\documentclass[10pt]{amsart}
\usepackage{graphicx}
\usepackage{url}

\usepackage{mathrsfs}

\newtheorem{theorem}{Theorem}[section]
\newtheorem{lemma}[theorem]{Lemma}           
\newtheorem{cor}[theorem]{Corollary}
\newtheorem{prop}[theorem]{Proposition}

\theoremstyle{definition}

\theoremstyle{remark}

\numberwithin{equation}{section}


\title[ITE problems on compact manifolds]
{Interior transmission eigenvalue problems on compact manifolds with  boundary conductivity parameters}
\author[H. Morioka]{Hisashi MORIOKA}
\address[H. Morioka]{Faculty of Science and Engineering,
Doshisha University, Tataramiyakodani 1-3, Kyotanabe, Kyoto, 610-0394, Japan}
\email{hmorioka@mail.doshisha.ac.jp}
\thanks{This work is partially supported by the JSPS grants-in-aid No. 16K17630 and by the Research Institute for Mathematical Sciences, a Joint Usage/Research Center located in Kyoto University.}

\author[N. Shoji]{Naotaka SHOJI}
\address[N. Shoji]{Kaichi Mirai Junior and Senior High School, Mugikura 1238, Kazo, Saitama, 349-1212, Japan
}
\email{nnao1003@math.tsukuba.ac.jp}

\date{\today}

\begin{document}
\maketitle

\begin{abstract} 
In this paper, we consider an interior transmission eigenvalue (ITE) problem on some compact $C^{\infty }$-Riemannian manifolds with a common smooth boundary.
In particular, these manifolds may have different topologies, but we impose some conditions of Riemannian metrics, indices of refraction and boundary conductivity parameters on the boundary. 
Then we prove the discreteness of the set of ITEs, the existence of infinitely many ITEs, and its Weyl type lower bound. 
For our settings, we can adopt the argument by Lakshtanov and Vainberg \cite{LaVa}, considering the Dirichlet-to-Neumann map.
As an application, we derive the existence of non-scattering energies for time-harmonic acoustic equations.
For the sake of simplicity, we consider the scattering theory on the Euclidean space. 
However, the argument is applicable for certain kinds of non-compact manifolds with ends on which we can define the scattering matrix.
\end{abstract}

%
%
\section{Introduction}

\subsection{Settings of ITE problems on manifolds}

We consider two connected and compact $C^{\infty} $-Riemannian manifolds $( M_1 , g_1 ) $ and $( M_2 , g_2 )$ with $C^{\infty }$-boundaries $\partial M_1 $ and $\partial M_2 $, respectively. 
We assume $d := \mathrm{dim} M_1 = \mathrm{dim} M_2  \geq 2 $ and $\mathrm{dim} \partial M_1 = \mathrm{dim} \partial M_2  = d-1$. 
Throughout of the present paper, we assume that 

\medskip

{\bf (A-1)} $M_1 $ and $M_2 $ have a common boundary $\Gamma := \partial M_1 = \partial M_2 $. 
$\Gamma $ is a disjoint union of a finite number of connected and closed components. 
The metrics satisfy $g_1 =g_2 $ on $\Gamma $.

\medskip

We will add some assumptions for $g_1 $ and $g_2 $ in a neighborhood of the boundary in \S 2.3. 
Note that we need our geometric assumptions only in some small neighborhoods of the boundary.
In particular, we do not assume that $M_1 $ and $M_2 $ are diffeomorphic outside of a small neighborhood of the boundary.

Let $\Delta_{g_k} $, $k=1,2 $, be the (negative) Laplace-Beltrami operator on each $M_k $.
We consider the following interior transmission eigenvalue (ITE) problem : 
\begin{gather}
 (-\Delta_{g_1}-\lambda n_1)u_1 = 0 \quad \text{in} \quad M_1, \label{S1_eq_ITE1}\\
 (-\Delta_{g_2}-\lambda n_2)u_2 = 0 \quad \text{in} \quad M_2, \label{S1_eq_ITE2}\\
 u_1-u_2 = 0, \quad \partial_{\nu_1}u_1-\partial_{\nu_2}u_2 = \zeta u_1 \quad \text{on} \quad \Gamma ,\label{S1_eq_ITE3} 
\end{gather} 
where each $n_k \in C^{\infty}(\overline{M_k})$, $k=1,2 $, is strictly positive on $M_k $ and $\zeta \in C^{\infty}(\Gamma)$. 
For $\zeta$, this paper covers the following two cases : (i) $\zeta =0$ on $\Gamma$, or (ii) $\zeta$ is strictly positive or strictly negative on every connected component of $\Gamma$.
Note that we also need a stronger assumption for $n_1 $, $n_2 $ or $\zeta $ in \S 3.2.

We call $\sqrt{n_k}$ and $\zeta$ the index of refraction on $M_k$ and the boundary conductivity parameter on $\Gamma$, respectively.
If there exists a pair of non-trivial solutions $(u_1 , u_2 ) \in H^2 (M_1 ) \times  H^2 ( M_2 )$ of (\ref{S1_eq_ITE1})-(\ref{S1_eq_ITE3}), we call corresponding $\lambda \in {\bf C} $ an \textit{interior transmission eigenvalue}.

\subsection{Backgrounds}
ITE problems naturally appears in inverse scattering problems for acoustic wave equations on ${\bf R}^d $ with compactly supported inhomogeneity.
In ${\bf R}^d $ for $ d\geq 2 $, time harmonic acoustic waves satisfy the equation 
\begin{equation}
(-\Delta -\lambda n )u=0 , \quad \lambda > 0 ,
\label{S1_eq_scatteringeq}
\end{equation}
where $n \in L^{\infty } ({\bf R}^d )$ is strictly positive in a bounded domain $\Omega $ with a suitable smooth boundary, and $n \big| _{{\bf R}^d \setminus \Omega } =1 $. 
Given an incident wave $u^i(x) = e^{i\sqrt{\lambda}x\cdot\omega}$ with an incident direction $\omega \in S^{d-1}$ and energy $\lambda > 0$, the scattered wave $u^s$ is described by the difference
between the total wave $u$ and the incident wave $u^i$ where $u$ is the solution of (\ref{S1_eq_scatteringeq}) satisfying the following asymptotic relation : as $ |x | \to \infty $ 
\begin{equation}
u(x) \simeq e^{i \sqrt{\lambda} x \cdot \omega } + C(\lambda ) |x|^{-(d-1)/2 } e^{ i \sqrt{ \lambda} |x| }  A( \lambda ; \omega , \theta ) , \quad \theta = x/|x|. 
\label{S1_eq_asymptoticexp}
\end{equation}
Here the second term on the right-hand side is the spherical wave scattered to the direction $\theta $. 
The function $A(\lambda ;  \omega , \theta )$ is the scattering amplitude.
The \textit{S-matrix} is given by $S(\lambda )= 1- 2\pi i  A(\lambda )$ where $A( \lambda )$ is an integral operator with the kernel $A(\lambda ; \omega , \theta )$.
Then the S-matrix is unitary operator on $L^2 (S ^{d-1 } )$.
If there exists a non zero function $\phi \in L^2 (S^{d-1} )$ such that $S(\lambda )\phi = \phi $ i.e. $A(\lambda ) \phi =0 $, we call $\lambda >0$ a \textit{non-scattering energy (NSE)}.
If $\lambda >0 $ is a NSE, we have that $u-u^i $ vanishes outside $\Omega $ from the Rellich type uniqueness theorem (see \cite{Re43} and \cite{Vek43}).
Hence we can reduce to the ITE problem 
\begin{gather}
( -\Delta -\lambda n ) v=0  \quad \text{in} \quad \Omega  , \label{S1_eq_ITE1euq} \\
( -\Delta -\lambda  ) w =0  \quad \text{in} \quad \Omega , \label{S1_eq_ITE2euq} \\
v = w , \quad  \partial_{\nu} v = \partial_{\nu} w \quad \text{on} \quad \partial \Omega , \label{S1_eq_ITEeuq} 
\end{gather} 
with $v=u$ and $w=u^i $. 
If $\lambda>0 $ is a NSE, $\lambda $ is also an ITE of the system (\ref{S1_eq_ITE1euq})-(\ref{S1_eq_ITEeuq}).
ITE problems were introduced in \cite{Ki} and \cite{CoMo} in the above view point.
For the Schr\"{o}dinger equation $(-\Delta +V -\lambda )u=0 $ with a compactly supported potential $V $ which satisfies $V(x) \geq \delta >0 $ in $\mathrm{supp} V$, we can state the ITE problem similarly.  
Recently, the ITE problem is generalized by \cite{Ves} to unbounded domains with exponentially decreasing perturbations at infinity.

As far as the authors know, results on the NSE are very scarce.
In particular, it seems to be no result for the existence of non-scattering energies except for spherically symmetric inhomogeneities (see \cite{CoMo}).
There are some examples of perturbations which do not have non-scattering energies (\cite{GeHa}, \cite{BlPaSy}, \cite{ElGu}, \cite{PaSaVe}). 
If the perturbation is compactly supported and the associated ITEs are discrete, the discreteness of NSE is a direct consequence.

The system (\ref{S1_eq_ITE1euq})-(\ref{S1_eq_ITEeuq}) is a kind of non self-adjoint problem. 
Moreover, we can construct a bilinear form associated with this system, but generally this bilinear form is not coercive. 
Note that the $T$-coercivity approach is valid for some anisotropic cases i.e. $-\Delta $ is replaced by $-\nabla \cdot A \nabla $ where $A$ is a strictly positive symmetric matrix valued function and $A \neq Id$. 
For the $T$-coercivity approach on this case, see \cite{BeChHa}.
Another common approach is to reduce an ITE problem to an equivalent forth-order equation.
For (\ref{S1_eq_ITE1euq})-(\ref{S1_eq_ITEeuq}), we can reduce to 
\begin{equation}
(\Delta + \lambda n) \frac{1}{n-1} (\Delta + \lambda ) \psi =0 ,\quad \psi = w-v \in H_0^2 (\Omega ) ,
\label{S1_eq_forthorder}
\end{equation}
which is formulated as the variational form 
\begin{equation}
 \int_{\Omega} \frac{1}{n-1} (\Delta\psi+\lambda\psi)(\Delta\overline{\phi}+\lambda n\overline{\phi})dx =0,
\label{S1_eq_variationforth}
\end{equation}
for any $ \phi \in H_0^2 (\Omega )$.
There are also many works on this approach for acoustic wave equations and Schr\"{o}dinger equations.
For more history, technical information and references on ITE problems, we recommend the survey by Cakoni and Haddar \cite{CaHa}.

\begin{figure}[t]
\centering
\includegraphics[width=10cm, bb=0 0 455 228]{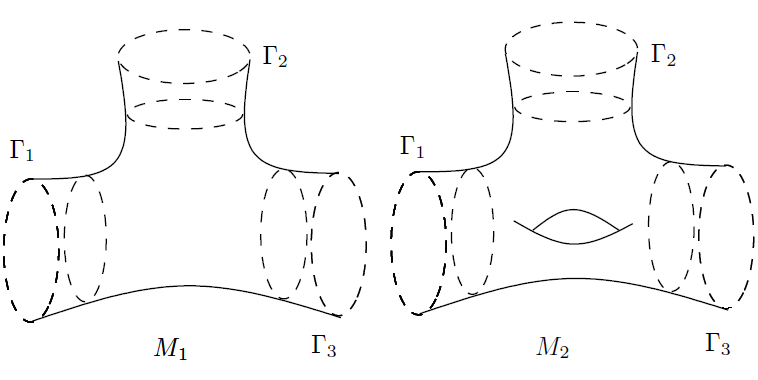}
\caption{Examples of $M_1 $ and $M_2 $ with common boundary $\Gamma = \bigcup _{j=1}^3 \Gamma_j $.}
\label{S1_fig_manifold}
\end{figure}

This paper consists of two parts. 
In the first part, we generalize the ITE problem in two directions.
The boundary conductivity parameter is introduced.
Moreover, we allow $M_1 $ and $M_2 $ to have different topologies (see Figure \ref{S1_fig_manifold}). 
We will discuss about ITEs in \S 2 and \S 3.

Forward and inverse scattering problems on non-compact manifolds are also well-known.
In particular, see e.g. \cite{Is} and \cite{IsKu} for asymptotically hyperbolic manifolds and see e.g. \cite{Ly}, \cite{Ly2}, \cite{IKL} and \cite{IKL2} for asymptotically cylindrical wavegudes.
We also mention that abundant references on related works are given in these articles.
Recently, the scattering theory on manifolds is derived by \cite{IS} without any assumptions on asymptotic behaviors of metrics.
We can define non-scattering energies on manifolds by the same way of the Euclidean space. 
Then the associated ITE problem on a compact manifold with a boundary is derived from the scattering theory on every manifold.
In particular, if we consider the scattering theory on a manifold with multiple ends, the associated bounded domain has multiple components of the boundary.

Since we do not assume that $ M_1 $ and $M_2 $ are diffeomorphic, it is difficult to use the forth-order equation approach.
Moreover, in view of assumptions (A-1) and (A-2) which is added in \S 2.3, the ITE problem is not elliptic, and we can not construct a suitable isomorphism $T$ such that the system (\ref{S1_eq_ITE1})-(\ref{S1_eq_ITE3}) is $T$-coercive. 
Therefore, neither the variational formulation approach nor the $T$-coercivity approach are valid for the proof of discreteness of ITEs in our case. 
Then we adopt arguments by Lakshtanov and Vainberg \cite{LaVa} in the present paper.
The approach in \cite{LaVa} is based on methods of elliptic pseudo-differential operators on the boundary and its application to the Dirichlet-to-Neumann (D-N) map.
The system (\ref{S1_eq_ITE1euq})-(\ref{S1_eq_ITEeuq}) is considered in \cite{LaVa}, but their argument is applicable to (\ref{S1_eq_ITE1})-(\ref{S1_eq_ITE3}) with the boundary conductivity parameter.
For the sake of the pseudo-differential calculus, we have imposed regularity conditions for $n_k$ and $\zeta $.

We should also mention about \cite{Vo} and \cite{PeVo}.
Recently, they proved the Weyl's asymptotics including complex ITEs and evaluated ITE-free regions in the complex plane under various conditions. 
They used the semi-classical analysis for the D-N map associated with an operator of the form $- n(x)^{-1} \nabla \cdot c(x) \nabla $ where $n,c$ are smooth and positive valued function on a bounded domain $\Omega \subset {\bf R}^d$.

In this paper, we construct the Poisson operator and the associated D-N map as elliptic pseudo-differential operators and we can compute exactly their symbols.
Using the ellipticity of the D-N map and the analytic Fredholm theory, we can prove the discreteness of the set of ITEs.
We also consider a Weyl type lower bound of the number of positive ITEs except for a small neighborhood of the origin.

A case which we can use the $T$-coercivity approach will be studied in the forthcoming paper \cite{MoSh}.

In the second part which will be discussed in \S 4, we derive the existence and a Weyl type lower bound of NSEs for the S-matrix of time-harmonic acoustic equations with compactly supported inhomogeneities.
In this paper, we consider the scattering theory on the Euclidean space for the sake of simplicity.
However, our argument is applicable to some kind of non-compact manifolds with ends (for example, Euclidean or hyperbolic ends) on which we can derive the scattering theory for suitable self-adjoint operators.
The main instrument is the equivalence of the S-matrix and the D-N map where the D-N map is defined for the interior Dirichlet problem in the support of the inhomogeneity.
This fact is often used in order to reduce inverse scattering problems (ISP) to inverse boundary value problems (IBVP). 
For this topic, see e.g. \cite{IsNa}, \cite{Is}, \cite{IsKu}, \cite{Es} and references therein.
Similarly, we reduce NSEs to ITEs.
In studies of ISP and IBVP, we can usually avoid Dirichlet eigenvalues associated with the interior Dirichlet problem.
However, we have to consider the Dirichlet eigenvalues for the study of NSEs and ITEs.
Hence we need to modify the proof of the equivalence of the S-matrix and the D-N map.

\subsection{Plan of the paper}

The plan of the paper is as follows.
In \S 2, we recall some basic properties of the D-N map.
For our purpose, we need to study about residues and regular parts of the D-N map near its poles. 
The relation between ITEs and non-trivial kernels of the difference of D-N maps is also introduced here.
Finally, we construct an approximate solution of the Dirichlet boundary value problems as a pseudo-differential operator, and we compute the symbol of the D-N map.
We prove our main results in \S 3. 
We use the analytic Fredholm theory, the parameter ellipticity of pseudo-differential operators and Weyl type asymptotic estimates for the number of Dirichlet eigenvalues on compact manifolds.
Our main results are Theorem \ref{S3_thm_discreteness} for the discreteness of ITEs and Theorem \ref{S3_thm_WeylNT} for the lower bound of the number of ITEs in $(\alpha , \infty)$ with sufficiently small $ \alpha >0 $. 
We discuss NSEs in \S 4.
After recalling some basic materials of the scattering theory, we prove the equivalence of the S-matrix and the D-N map, considering exterior and interior Dirichlet problems.

\subsection{Notation}
We use the following notations. 
We put ${\bf R} _{\geq 0} := [0,\infty ) $ and ${\bf R} _{>0} := (0, \infty ) $.
For the Riemannian metric $g _k = (g_{k,ij} ) $ of $M_k $, $\sqrt{g_k} $ and $ (g^{ij}_k )$ denote $\sqrt{\mathrm{det} g_k}$ and $g_k^{-1} $, respectively.
$dV_k (x)  : = \sqrt{g_k} dx $ and $dS  (x) $ denote the volume element on $M_k $ and the surface element on $\Gamma $ induced by $dV_k (x) $, respectively.
We often write them as $ dV_k $ and $dS$ omitting $(x)$.
Letting $x= (x_1 , \cdots , x_d )$ be a local coordinate of $M_k$, $\partial_j $ or $\partial _{x_j} $ denote $\partial / \partial x_j $.
For $\xi $, we use the similar manner.
For a multiple index $ \alpha = ( \alpha_1 , \cdots , \alpha _d ) $, we write $ \partial ^{\alpha }_x = \partial _1^{\alpha_1} \cdots \partial_d ^{\alpha_d } $. 
We often compute some kind of symbols $p(x,\xi )$.
For short, we denote by $p(x, i \partial_x )$ a pseudo-differential operator where each $\xi _j $ of $p(x,\xi )$ is replaced by $ i\partial _{x_j} $. 
Similarly, when we write $p(-i\partial _{\xi} , \xi )$, each $x_j$ is replaced by $-i\partial _{\xi_j} $.
$ \partial _{\nu_k} $ denotes the outward normal derivative on $\Gamma $ associated with $M_k$.
For a strictly positive valued function $\eta \in L^{\infty} (M_k) $, $L^2 ( M_k , \eta dV_k )$ is the $L^2 $ space on $M_k $ with the inner product $(u,v ) _{L^2 (M_k , \eta dV_k )} = ( \eta u ,v) _{L^2 (M_k )} $.


\section{Dirichlet-to-Neumann map}
 \subsection{Dirichlet-to-Neumann map}
Here we consider the following Dirichlet problems :
\begin{equation}
( -\Delta_{g_k} -\lambda n_k  ) u_k =0  \quad \text{in} \quad M_k , \quad u_k = f \quad \text{on} \quad \Gamma  , \label{S1_eq_D1} 
\end{equation} 
for $k=1,2 $.
We define the Dirichlet-to-Neumann (D-N) map $\Lambda _k (\lambda ) $ by 
\begin{equation}
\Lambda _k (\lambda ) f = \partial_{\nu_k } u_k \quad \text{on} \quad \Gamma ,
\label{S2_def_DN}
\end{equation}
where $u_k $ is a solution of (\ref{S1_eq_D1}).

In the following, we call $\lambda  $ a Dirichlet eigenvalue if there exists a non-trivial solution of the equation 
\begin{equation}
( - \Delta _{g_k} -\lambda  n_k ) u_k =0  \quad \text{in} \quad M_k , \quad u_k =0 \quad \text{on} \quad \Gamma .
\label{S2_eq_dirichleteigen}
\end{equation}
In fact, (\ref{S2_eq_dirichleteigen}) is equivalent to 
\begin{equation}
( - n_k ^{-1} \Delta _{g_k} -\lambda  ) u_k =0 \quad \text{in} \quad M_k , \quad u_k =0 \quad \text{on} \quad \Gamma ,
\label{S2_eq_dirichleteigen2}
\end{equation}
which is an eigenvalue problem of the second-order self-adjoint elliptic operator $L_k = - n_k^{-1} \Delta _{g_k} $ in $L^2 (M_k , n_k dV_k )$ with the Dirichlet boundary condition.
Then its eigenvalues form an increasing sequence $ 0< \lambda_{k,1} \leq \lambda_{k,2} \leq \cdots $, satisfying the Weyl's asymptotics which we derive in \S 3.
The corresponding eigenfunctions $ \phi _{k,j} $ can be chosen so that $\{ \phi_{k,j} \} $ is an orthonormal basis in $L^2 (M_k , n_k dV_k  )$.
We denote the set of Dirichlet eigenvalues by $\{ \lambda _{k,j} \} := \{ \lambda _{k,j} \}_{j=1}^{\infty}   $.
For $\lambda \not\in \{ \lambda_{k,j} \} $, the D-N map $\Lambda _k (\lambda )$ is well-defined and extends uniquely as a continuous operator $\Lambda _k (\lambda ) : H^{3/2} (\Gamma ) \to H^{1/2} (\Gamma )$.

Let $ \mathcal{E}_{k,j} \subset {\bf Z}_+ $ such that $ \bigcup _{j=1}^{\infty} \mathcal{E}_{k,j} = {\bf Z}_+ $, and $ i_1 $ and $i_2 $ belong to the same set $ \mathcal{E}_{k,j} $ if and only if $ \lambda _{k,i_1} = \lambda _{k,i_2} $.
We denote eigenvalues corresponding $ \mathcal{E}_{k,j} $ by $\lambda _{k,(j)} $.
$\mathcal{L}(\lambda_{k,i})$ means the set $\mathcal{E}_{k,j}$ with $\lambda_{k,(j)} = \lambda_{k,i}$

\begin{prop}
$ \Lambda_k ( \lambda ) $ is meromorphic with respect to $\lambda \in {\bf C} $ and has first order poles at $\lambda \in \{ \lambda_{k,j} \} $.
Moreover, $\Lambda _k (\lambda )$ has the following representations : \\
(1) For $x\in \Gamma $ and $f\in H^{3/2} (\Gamma )$, we have
 \begin{equation}
 \Lambda _k (\lambda )f  (x)
= -\int_{\Gamma} \sum_{j=1}^{\infty} \frac{  \partial _{\nu _k (x)}  \phi_{k,j} (x) \,   \partial_{\nu _k (y)} \phi_{k,j} (y)  }{\lambda_{k,j} -\lambda } f(y) dS (y).
 \label{S2_eq_DNrepresentation}
 \end{equation}
 (2) In a neighborhood of $\lambda_{k,j}$, we have
 \begin{equation}
\Lambda _k (\lambda )= \frac{Q_{k, \mathcal{L}(\lambda_{k,j} )}}{\lambda_{k,j} -\lambda } + H_k (\lambda ) ,
\label{S2_eq_resDN}
\end{equation}
where $Q_{k, \mathcal{L}( \lambda_{k,j} ) } $ is the residue of $\Lambda_k (\lambda )$ at $ \lambda = \lambda_{k,j} $ given by 
\begin{equation}
Q_{k, \mathcal{L}( \lambda_{k,j} ) } f=  - \sum_{i\in \mathcal{L} (\lambda_{k,j} )}  \int _{\Gamma} \partial_{\nu_k (y)} \phi _{k,i } (y) \,  f (y)  dS (y) \, \partial_{\nu_k } \phi _{k,i } ,
\label{S2_eq_Pk}
 \end{equation}
 and $H_k (\lambda ) : H^{3/2} (\Gamma ) \to H^{1/2} (\Gamma ) $ is analytic in a neighborhood of $\lambda_{k,j} $.

\label{S2_prop_green_meromor}
\end{prop}

Proof. 
We can follow the argument of \S 4.1.12 in \cite{KaKuLa}.
Let $E_k \in H^2(M_k)$ be an extension of $f$ into $M_k$ satisfying $E_k \big|_{\Gamma} = f$ and $\|E_k\|_{H^2(M_k)} \le C\|f\|_{H^{3/2}(\Gamma)}$ for some constants $C > 0$. 
Then we have
$$
 (-n_k^{-1}\Delta_{g_k}-\lambda)(u_k-E_k) = -(-n_k^{-1}\Delta_{g_k}-\lambda)E_k ,
$$
where $u_k$ is a solution of (\ref{S1_eq_D1}).
Since $R_k(\lambda) := (-n_k^{-1}\Delta_{g_k}-\lambda)^{-1}$ is a meromorphic operator valued function with first order poles only at $\lambda \in \{\lambda_{k,j}\}$, $u_k = E_k-R_k(\lambda)(-n_k^{-1}\Delta_{g_k}-\lambda)E_k$ is also a meromorphic $H^2(M_k)$-valued function with first order poles only at $\lambda \in \{\lambda_{k,j}\}$.

Next we prove (\ref{S2_eq_DNrepresentation}).
Integrating by parts, we compute the Fourier coefficients of $u_k$ with respect to the real-valued eigenfunction $\phi_{k,j}$ :
\begin{equation}
 (u_k,\phi_{k,j})_{L^2(M_k,n_kdV_k)} = -\int_{\Gamma} \frac{ \partial_{\nu_k(y)}\phi_{k,j}(y)}{\lambda_{k,j}-\lambda} f(y) dS(y).
 \label{S2_eq_fouriercoefficients}
\end{equation}
From this formula and the outward normal derivative of $u_k$, $\Lambda_k(\lambda)$ satisfies (\ref{S2_eq_DNrepresentation}).

Finally we verify (\ref{S2_eq_resDN}) and (\ref{S2_eq_Pk}).
Let $P_{k,j} : L^2 (M_k , n_k dV_k ) \to L^2 (M_k , n_k dV_k ) $ be the projection to the eigenspace corresponding to $\lambda _{k,j}$ i.e.   
$$
P_{k,j} v = \sum_{i\in \mathcal{L} (\lambda_{k,j} ) } ( v , \phi _{k,i} ) _{L^2 (M_k , n_k dV_k ) } \phi _{k,i} , \quad v\in L^2 (M_k , n_k dV_k )  .
$$ 
In view of (\ref{S2_eq_fouriercoefficients}), we have 
$$
P_{k,j} u_k  = -\frac{1}{\lambda_{k,j} -\lambda }\sum_{i \in \mathcal{L} (\lambda_{k,j} )}\int _{\Gamma} \partial_{\nu_k  (y)} \phi_{k,i} (y) \, f(y) dS (y) \, \phi _{k,i }  ,
$$ 
and this implies (\ref{S2_eq_Pk}). 
Moreover, 
$$
 (1-P_{k,j} ) u_k = -\sum_{i\not\in \mathcal{L} (\lambda_{k,j} )} \frac{1}{\lambda_{k,i} -\lambda }  \int _{\Gamma}  \partial_{\nu _k (y)} \phi_{k,i} (y) \, f(y) dS (y) \, \phi _{k,i }  ,
  $$
 is analytic with respect to $\lambda $ in a neighborhood of $\lambda _{k,j} $. 
Putting $ H_k (\lambda )f = \partial _{\nu_k} ( (1-P_{k,j} ) u_k ) $ on $\Gamma $, we have the proposition.
\qed

\medskip

\textit{Remark.} 
The formula (\ref{S2_eq_Pk}) means that the range of $Q_{k, \mathcal{L}( \lambda_{k,j} ) } $ is a finite dimensional subspace spanned by $ \partial _{\nu _k} \phi _{k,i} $ for $i \in \mathcal{L} (\lambda _{k,j} )$.
Note that $\partial _{\nu_k } \phi _{k,i} $ for all $i\in \mathcal{L} (\lambda _{k,j} ) $ are linear independent since $\phi _{k,i} $ are orthogonal basis.
Hence $\mathrm{dim} \mathrm{Ran} Q_{k,\mathcal{L} (\lambda _{k,j} )} $ coincides with the multiplicity of $\lambda _{k,j} $. 
We can see that the integral kernel of $Q_{k, \mathcal{L}( \lambda_{k,j} ) } $ is smooth in $(x,y)$ by the regularity property of Dirichlet eigenfunctions.

\medskip

As has been in Propositions \ref{S2_prop_green_meromor}, $\Lambda_1 (\lambda ) - \Lambda_2 (\lambda )$ is also meromorphic with respect to $\lambda \in {\bf C}$ and has first order poles at $\lambda \in \{\lambda_{1,j}\}\cup\{\lambda_{2,j}\}$. 
In a neighborhood of a pole $\lambda_0 $, we have 
\begin{equation}
\Lambda_1 (\lambda ) - \Lambda_2 (\lambda )= \frac{Q _{\lambda_0} }{\lambda_0 -\lambda }+ H_{\lambda_0 } (\lambda ) ,
\label{S2_eq_laurentdiff}
\end{equation}
where $Q _{\lambda_0}  $ and $H_{\lambda_0 } (\lambda ) $ have same properties of $Q_{k, \mathcal{L}( \lambda_{k,j} ) }  $ and $H_k (\lambda )$, respectively. 
In the following, we define the kernel of $\Lambda_1 ( \lambda ) - \Lambda _2 (\lambda )$ by 
\begin{gather}
\begin{split}
&\mathrm{Ker} (\Lambda_1 (\lambda )-\Lambda_2 (\lambda ) ) \\
&= \left\{ 
\begin{split}
& \{ f\in H^{3/2} (\Gamma ) \ ; \  (\Lambda_1 (\lambda )-\Lambda_2 (\lambda ) ) f=0 \} , \ \text{if} \ \lambda \ \text{is not a pole}, \\
& \{ f\in H^{3/2} (\Gamma ) \ ; \  Q_{\lambda_0 } f = H_{\lambda_0} (\lambda_0 )f =0 \} , \ \text{if} \ \lambda =\lambda_0 \ \text{is a pole}.
\end{split}
\right.
\end{split}
\label{S2_def_kernelDNdiff}
\end{gather}
For $ \Lambda_1 (\lambda ) -\Lambda_2 (\lambda ) -\zeta $, we define its kernel by the same manner.

\begin{lemma}
Let $\lambda \in \{ \lambda_{k,j} \} $.
Then the equation (\ref{S1_eq_D1}) has a non trivial solution if and only if $f $ is orthogonal to $\partial _{\nu_k } \phi _{k,j} $ in $L^2 (\Gamma )$ for all $j \in \mathcal{L} (\lambda )$. 
\label{S2_lem_singular_sol}
\end{lemma}

Proof. 
If $f $ is orthogonal to $\partial _{\nu_k } \phi _{k,j} $ for all $j\in \mathcal{L} (\lambda )$, there exist general solutions of the form 
\begin{equation}
u_k = -\sum_{i\not\in \mathcal{L} (\lambda ) } \frac{1}{\lambda_{k,i} -\lambda } \int_{\Gamma} \partial _{\nu_k (y)} \phi _{k,i} (y) \, f(y) dS (y) \,  \phi _{k,i} + \sum_{i\in \mathcal{L} (\lambda ) } c_i \phi _{k,i} ,
\label{S2_eq_singularsol2}
\end{equation}
for any $c_i \in {\bf C} $. 

If $u_k $ is a non trivial solution of (\ref{S1_eq_D1}), we have by Green's formula 
$$
\int _{M_k} \left(  \Delta_{g_k} u_k  \cdot \phi _{k,i} - u_k \cdot \Delta_{g_k} \phi_{k,i} \right) dV_k  = - \int_{\Gamma}  u_k \cdot \partial_{\nu_k } \phi _{k,i} dS, 
$$
for $i\in \mathcal{L} (\lambda )$.
Since $\lambda = \lambda _{k,i} $, the left-hand side is equal to zero.
Then $f = u_k \big| _{\Gamma} $ is orthogonal to $\partial _{\nu_k} \phi _{k,i} $.
\qed

\medskip

The above lemma implies a unique solvability in a subspace as follows.

\begin{cor}

Let $E_k (\lambda_0 ) \subset H^2 (M_k )$ be the eigenspace spanned by $ \phi _{k,i} $, and $B_k (\lambda_0 )$ be the subspace of $H^{3/2} (\Gamma )$ spanned by $\partial _{\nu_k} \phi _{k,i} $ for all $i\in \mathcal{L} (\lambda_0 )$ with $\lambda_0 \in \{ \lambda_{k,j} \}$. 
We denote by $E_k (\lambda_0 )^c $ and $B_k (\lambda_0 )^c $ their orthogonal complements in $L^2 (M_k)  $ and $L^2 (\Gamma )$, respectively. 
For any $f\in B_k (\lambda_0 )^c$, there  exists a unique solution $u_k \in E_k (\lambda_0 ) ^c \cap H^2 (M_k ) $ of (\ref{S1_eq_D1}) represented by 
\begin{equation}
u_k = -\sum_{i\not\in \mathcal{L} (\lambda _0 ) } \frac{1}{\lambda_{k,i} -\lambda } \int_{\Gamma}  \partial _{\nu_k (y)} \phi _{k,i} (y) \,  f(y) dS (y) \,  \phi _{k,i} .
\label{S2_eq_singularsol11}
\end{equation}
\label{S2_cor_singular_sol2}
\end{cor}

Proof.
We have only to check the uniqueness.
This is trivial since the equation (\ref{S2_eq_dirichleteigen}) has only the trivial solution in $E _k (\lambda_0 )^c $. 
\qed

\medskip

Now we can state the relation between ITEs and the D-N map as follows.

\begin{lemma}
(1) Suppose $\lambda \not\in \{ \lambda_{1,j} \} \cap \{ \lambda_{2,j} \} $.
Then $\lambda \in {\bf C} $ is an ITE if and only if $ \mathrm{Ker} (\Lambda _1 (\lambda ) - \Lambda _2 (\lambda ) -\zeta ) \not= \{ 0 \}$.
The multiplicity of $\lambda $ coincides with $\mathrm{dim} ( \mathrm{Ker} (\Lambda _1 (\lambda ) - \Lambda _2 (\lambda ) -\zeta ) )$. \\
(2) Suppose $\lambda \in \{ \lambda_{1,j} \} \cap \{ \lambda_{2,j} \} $.
Then $\lambda $ is an ITE if and only if $ \mathrm{Ker} (\Lambda _1 (\lambda ) - \Lambda _2 (\lambda )- \zeta ) \not= \{ 0 \} $ or the ranges of $ Q_{1,\mathcal{L} (\lambda )} $ and $Q_{2, \mathcal{L} (\lambda )} $ have a non trivial intersection. 
The multiplicity of $\lambda $ coincides with the sum of $\mathrm{dim} ( \mathrm{Ker} (\Lambda _1 (\lambda ) - \Lambda _2 (\lambda )- \zeta ) )$ and the dimension of the above intersection.  
\label{S2_lem_ITEeqcondition}
\end{lemma}

Proof. 
We first prove the assertion (1).
When $\lambda \not\in \{ \lambda_{1,j} \} \cup \{ \lambda_{2,j} \} $, this lemma is a direct consequence of the definition of ITEs.
We have only to show for $\lambda \in \{ \lambda_{1,j} \} \setminus \{ \lambda_{2,j} \}$. 
For $ 0 \not= f \in \mathrm{Ker} ( \Lambda_1 (\lambda ) - \Lambda_2 (\lambda )-\zeta ) $, we have $ Q_{1, \mathcal{L} (\lambda )} f = (H_1 (\lambda ) - \Lambda_2 (\lambda )-\zeta )f=0 $.
From $Q_{1, \mathcal{L} (\lambda )} f = 0$ and (\ref{S2_eq_Pk}), we have $f \in B_1 (\lambda )^c $. 
By Lemma \ref{S2_lem_singular_sol} and Corollary \ref{S2_cor_singular_sol2}, the following equation has a unique non trivial solution : 
\begin{gather}
 (-\Delta_{g_1}-\lambda n_1)u_1 = 0 \quad \text{in} \quad M_1, \quad u_1 = f \quad \text{on} \quad \Gamma.
 \label{S2_eq_DD1}
\end{gather}
On the other hand, from $(H_1(\lambda)-\Lambda_2(\lambda)-\zeta ) f=0$, we have 
\begin{gather}
 (-\Delta_{g_2}-\lambda n_2)u_2 = 0 \quad \text{in} \quad M_2, \quad u_2 = f,\ \partial_{\nu_2}u_2 = (H_1(\lambda)-\zeta ) f \quad \text{on} \quad \Gamma.
 \label{S2_eq_DD2}
\end{gather}
Summarizing (\ref{S2_eq_DD1}) and (\ref{S2_eq_DD2}) and $\partial_{\nu_1}u_1 = H_1(\lambda)f$, $\lambda$ is an ITE.
Conversely, if $\lambda $ is an ITE, from Lemma \ref{S2_lem_singular_sol}, the equation (\ref{S1_eq_D1}), $k=1 $, with the condition $u_1 \big| _{\Gamma} =f \not= 0 $ must have a non trivial solution.
In view of (\ref{S2_eq_Pk}), we have $f\in B_1 (\lambda )^c $, and this implies $Q_{1, \mathcal{L} (\lambda )} f=0 $.
This means $\partial_{\nu_1} u_1 =  H_1(\lambda) f$.
On the other hand, $\partial _{\nu _1 } u_1 - \partial_{\nu_2 } u_2 = \zeta f$ means $ (H_1 (\lambda ) - \Lambda_2 (\lambda )-\zeta )f=0$.
Therefore, $f$ must be in $\mathrm{Ker} (\Lambda_1 (\lambda ) - \Lambda_2 (\lambda )-\zeta )$.
We have proven the assertion (1).

For the assertion (2), we have only to show the latter case.
In fact, if there exists a non trivial solution $(u_1  , u_2 )$ of 
\begin{gather*}
(-\Delta_{g_1} -\lambda n_1 ) u_1 =0  \quad \text{in} \quad M_1 , \\
(-\Delta_{g_2} -\lambda n_2 ) u_2 =0  \quad \text{in} \quad M_2 ,
\end{gather*} 
with $u_1 = u_2   =0 $ and $\partial _{\nu_1} u_1 = \partial_{\nu_2 } u_2  $ on $\Gamma$, then we have that the ranges of $ Q_{1,\mathcal{L} (\lambda )} $ and $Q_{2, \mathcal{L} (\lambda )} $ have a non trivial intersection, recalling $\mathrm{Ran} Q_{k,\mathcal{L}(\lambda)} = \mathrm{Span} \{ \partial_{\nu_k} \phi_{k,j} \}_{j \in \mathcal{L}(\lambda)}$ for $k = 1,2$.  
Conversely, if the ranges of $ Q_{1,\mathcal{L} (\lambda )} $ and $Q_{2, \mathcal{L} (\lambda )} $ have a non trivial intersection, then there exists a non trivial solution $(u_1  , u_2 )$ of the above system with the condition $u_1= u_2   =0 $ and $\partial _{\nu_1} u_1 = \partial_{\nu _2} u_2  $ on $\Gamma$, since $\partial _{\nu_k} \phi _{k,i} $ for all $i\in \mathcal{L}(\lambda_0)$ are linear independent. 
Then $\lambda $ is an ITE. 
\qed

\medskip

\textit{Remark.} In \cite{LaVa}, the authors call $\lambda$ \textit{singular ITE} if $\lambda $ satisfies the latter condition in the assertion (2) of Lemma \ref{S2_lem_ITEeqcondition}.

%
%

\subsection{Parametrix}
Now let us compute the symbol of the D-N map. 
Here we construct the parametrix for (\ref{S1_eq_D1}).
As in \cite{LaVa}, we follow the argument of \S 2 in \cite{VaGu}, slightly modifying it for our case.

In the following, we assume that the equation (\ref{S1_eq_D1}) is uniquely solvable in $H^2 (M_k )$ or a suitable subspace of $L^2 (M_k )$.

We take a point $x^{(0)} \in \Gamma $ and fix it.
Let $V \subset \Gamma $ be a sufficiently small neighborhood of $x^{(0)} $ in $\Gamma $. 
There  exist small open domains $U_k \subset M_k $, $k =1,2$, such that $ \overline{U_k} \cap \Gamma = V$ and $U_1 $ and $U_2 $ are diffeomorphic to an open domain $\Omega \subset {\bf R}^d $.

We introduce local coordinates $y=(y_1 , \cdots , y_{d-1} , y_d ) $ in $\Omega$ with the center $x^{(0)} \in V$ such that $x^{(0)} = 0$, $\Omega $ is given by $y_d >0 $, $|y| < \epsilon _0 $ for a small $\epsilon_0 >0 $, the subset $\partial \Omega ^0 := \{ y\in \overline{\Omega} \ ; \ y_d =0 \} $ is diffeomorphic to $V$, and $y_d $ is the distance between a point $y=(y_1 , \cdots , y_{d-1}  ,y_d) \in \Omega $ and $\partial \Omega^0 $. 
Then $y=(y_1 ,\cdots ,y_d)$ are common local coordinates of $U_1 $ and $U_2$.
Therefore, we have 
$$
( g^{ij}_k (y) )_{i,j} = \left[ \begin{array}{cc}
\widetilde{g}_k (y') & \widetilde{p}_k (y) \\
{}^t \widetilde{p}_k (y) & 1 \end{array} \right] , \quad y' = (y_1 , \cdots , y_{d-1} ) ,
$$
in $U_k$ where $ \widetilde{g}_k (y' )= ( \widetilde{g} _k^{ij} (y') )_{i,j} $ is a smooth, positive definite and symmetric $(d-1)\times(d-1)$-matrix valued function, and $ \widetilde{p} _k (y)= {}^t ( p_{k,1} (y ) , \cdots , \widetilde{p} _{k,d-1} (y) )$ is a $(d-1)$-dimensional vector valued function.

A function $ F( y' , y_d , \xi ' , \xi _d ) $ with $(y' , y_d ) , ( \xi ' , \xi_d ) \in {\bf R}^{d} $ is homogeneous of the generalized degree $s $ if $F $ satisfies 
\begin{equation}
F( t ^{-1} y' , t ^{-1} y_d , t \xi ' , t \xi _d )= t^s F( y' , y_d , \xi ' , \xi _d ) , 
\label{S2_def_hom_gen}
\end{equation}
for any $t>0 $.
For $ F(y_d , \xi ') $, we define the homogeneity by the similar manner.

Taking the $y$-coordinate as above, we can rewrite $ A_k = -\Delta_{g_k} -\lambda n_k $ as 
\begin{equation}
A_k 
=: - \partial_d^2 - \sum_{i,j=1}^{d-1} \widetilde{g}_k^{ij} (y') \partial_i \partial_j  -2 \sum_{i=1}^d \widetilde{p} _{k,i} (y) \partial _i \partial_d - \sum_{i=1}^d \widetilde{h}_{k,i} (y  ) \partial_i -\lambda n_k (y) ,
\label{S2_def_Ak}
\end{equation}
in $U_k $ with $  \widetilde{h} _{k,i } (y) = (\sqrt{ g_k } )^{-1}  \sum_{j=1}^{d} \partial_j (\sqrt{g_k } \, g^{ij}_k ) $. 
Note that $ \widetilde{g}_k^{ij} (y') $, $ \widetilde{p} _{k,i} (y) $ and $ \widetilde{h} _{k,i} (y)$ are defined by $g _k (y)$.
In view of the assumption (A-1), we have in $y$-coordinates that $ \widetilde{g} _1^{ij} (y') = \widetilde{g}_2 ^{ij} (y') $, $\widetilde{p}_{1,i}(y) |_{y_d = 0} = \widetilde{p}_{2,i}(y) |_{y_d = 0} =0$.

The symbol of $A_k $ is given by 
\begin{gather}
\begin{split}
&A_k(\lambda; y' , y_d , \xi ' , \xi_d ) \\
&= \xi _d^2 +\sum_{i,j =1 }^{d-1} \widetilde{g}_k^{ij} (y') \xi_i \xi_j +2 \sum_{i=1}^d \widetilde{p}_{k,i} (y) \xi _i \xi_d -i \sum_{i=1}^{d-1} \widetilde{h}_{k,i}  (y ) \xi_i -\lambda n_k (y)  .
\end{split}
\label{S2_eq_symbolAk}
\end{gather}

In the following, let $N>0 $ be a sufficiently large integer.
Now we take $z = (z' ,0) \in \partial \Omega ^0 $ arbitrarily and fix it. 
Using the Taylor series of $ \widetilde{g}_k^{ij} (y') $, $\widetilde{p}_{k,i} (y)$, $\widetilde{h}_{k,i}  (y) $ and $n_k (y) $ with respect to $y$ centered at $(z' ,0 ) \in \partial \Omega ^0 $, we can expand the symbol $A _k (y' , y_d , \xi ' , \xi_d )$ of $A_k $ as the sum of following terms :
\begin{equation}
\xi_d^2 + \sum_{i,j =1}^{d-1} \widetilde{g} _k^{ij} (z')  \xi _i \xi _j , \label{S2_eq_expand_symbol1} 
\end{equation}
\begin{equation}
\begin{split}
& \sum_{i,j=1}^{d-1} \nabla _{y'} \widetilde{g}_k^{ij} (z') \cdot (y'-z') \xi _i \xi _j  + i \sum _{i=1}^d \widetilde{h}_{k,i} (z',0 ) \xi_i \\
&  +2 \sum_{i=1}^{d-1} \left( \nabla_{y'} \widetilde{p} _{k,i} (z' ,0) \cdot (y'-z' ) + \partial_d \widetilde{p}_{k,i} (z' ,0 ) y_d \right) \xi_i \xi_d , 
\end{split}
\label{S2_eq_expand_symbol2}
\end{equation}
and
\begin{equation}
\begin{split}
&\sum_{i,j=1}^{d-1} \sum_{|\alpha '| =m} \frac{\partial _{y'}^{\alpha '} \widetilde{g}_k^{ij} (z')}{\alpha ' !}  (y'-z') ^{\alpha '} \xi_i \xi_j + 2 \sum_{i=1}^d \sum_{ | \alpha | =m } \frac{\partial_y^{\alpha} \widetilde{p}_{k,i} (z',0)}{\alpha !}  (y'-z')^{\alpha '} y_d^{\alpha_d} \xi_i \xi _d  \\
&+i \sum_{i=1}^d \sum _{| \alpha | =m-1} \frac{ \partial _y^{\alpha} \widetilde{h} _{k,i} (z',0 )}{\alpha !} (y'-z')^{\alpha '} y_d^{\alpha_d} \xi_i -\lambda  \sum_{|\alpha | =m-2} \frac{ \partial _y^{\alpha} n_k (z',0 ) }{\alpha !} (y'-z') ^{\alpha '} y_d^{\alpha_d} ,
\end{split}
 \label{S2_eq_expand_symbol3}
\end{equation}
for $2 \leq m \leq N $ with the remainder term which has zero of order $N-1$ at $y' =0$ or $ (y' , y_d )= (0,0) $. 
We rewrite the sum of (\ref{S2_eq_expand_symbol1})-(\ref{S2_eq_expand_symbol3}) and the remainder term as 
\begin{gather}
\begin{split}
& A_k(\lambda;y',y_d,\xi',\xi_d) = A_{k,0} ( z' ; \xi ' , \xi _d ) + A_{k,1} (  z' ; y'-z' , y_d , \xi ' , \xi_d ) \\
&+ \sum_{m=2}^N A_{k,m} (\lambda , z' ; y' - z' , y_d , \xi '  , \xi_d  ) 
+ A'_{k,N+1} ( \lambda , z' ; y' - z' , y_d , \xi '  , \xi_d ) .
\end{split}
\label{S2_eq_expansion}
\end{gather}
Then each $A_{k,m} $ is a homogeneous polynomial in $y' -z' , y_d , \xi' , \xi_d $ of generalized degree $2-m $.
In particular, $ A_{k,0} $ is the principal symbol of $A_k $. 
$A' _{k,N+1} $ vanishes at $(z' ,0)$ and the order of the zero is $N-1$.

In the following arguments, we put
\begin{equation}
|\xi '|_{\Gamma} ^2 := \sum_{i,j =1 }^{d-1} \widetilde{g}_k^{ij} (y' ) \xi_i \xi_j .
\label{S2_def_gstar}
\end{equation}
We define the following differential operators :  
\begin{equation}
\widetilde{A}_{k,0} = A_{k,0} (   z' ;  \xi' , i \partial_d ) = - \partial_d^2 + |\xi '| _{\Gamma}^2 , 
\label{S2_def_Atilde}
\end{equation}
\begin{equation}
\widetilde{A}_{k,1} = A_{k,1} (   z' ;  -i \partial _{\xi '} , y_d , \xi ' , i \partial _d  )  , 
\label{S2_def_Atilde2}
\end{equation}
and
\begin{equation}
\widetilde{A}_{k,m} = A_{k,m} ( \lambda ,   z' ;  -i \partial _{\xi '} , y_d , \xi ' , i \partial _d ) , \quad m\geq 2 . 
\label{S2_def_Atilde3}
\end{equation}


\begin{prop}
Let $F(y_d , \xi ' )$ be a smooth function and homogeneous of the generalized degree $s$ with respect to $y_d $ and $\xi ' $.
Then we have that $\widetilde{A} _{k,m} F $ is the homogeneous of the generalized degree $2-m+s $ with respect to $y_d $ and $\xi ' $.
\label{S2_prop_degree}
\end{prop}

Proof.
Note that $F (y_d , \xi ' ) = |\xi '| ^s F ( |\xi '| y_d , \xi ' / | \xi '| )$.
Then we can show that $ \partial_d F$ and $ \partial _{\xi_j} F$ are homogeneous of generalized degree $s+1$ and $s-1 $, respectively.
\qed

\medskip

Now let us construct an approximate solution of (\ref{S1_eq_D1}).

\begin{lemma}
Suppose $ |\xi ' | _{\Gamma}  \not= 0$.
The system of second order ordinary differential equations
\begin{gather}
\widetilde{A}_{k,0} E_{k,0} (z' ; y_d , \xi ' )=0 , \label{S2_eq_sys0} \\
\widetilde{A}_{k,0} E_{k,1} (z' ; y_d , \xi ' )= -\widetilde{A}_{k,1} E_{k,0} (z' ; y_d , \xi' ) ,\label{S2_eq_sys1} \\
\cdots \nonumber \\
\widetilde{A}_{k,0} E_{k,m} (  z' ; y_d , \xi' )= -\sum_{n=1}^m \widetilde{A}_{k,n} E_{k,m-n} ( z' ; y_d , \xi ' )  , \label{S2_eq_sys2}
\end{gather}
has a unique solution $\{ E_{k,m} \} _{m=0,1,2,\cdots }$ such that each $E_{k,m} $ converges to zero as $y_d \to \infty $ and satisfies 
$$
E_{k,0} \big| _{y_d =0 } =1, \quad E_{k,m} \big| _{y_d =0 } =0, \quad m\geq 1 .
$$
In particular, we have $E_{k,0} (z' ; y_d , \xi ')= e^{-|\xi ' |_{\Gamma} y_d }$.
Each solution $E_{k,m}  $ is smooth and homogeneous with respect to $y_d $ and $\xi '$ of generalized degree $-m$. 
(For $m \geq 2 $, each $E_{k,m} $ depends also on $\lambda $. We omit $\lambda $ in the notation.)
\label{S2_lem_systemODE}
\end{lemma}

Proof.
Since $ \widetilde{A}_{k,0} = - \partial_d^2 +|\xi '|^2_{\Gamma} $, we have $E_{k,0} (z' ; y_d , \xi ')= e^{-|\xi ' |_{\Gamma} y_d }$.
Obviously, $E_{k,0} $ is homogeneous of the generalized degree $0$.
Let us consider the equation
\begin{equation}
( -\partial_d^2 + |\xi ' |^2_{\Gamma} ) v = p \quad \text{on} \quad (0, \infty ) ,
\label{S2_eq_ode_helmholtz}
\end{equation}
for $v(y_d, \xi ') $ and $ p (y_d , \xi ')$ with $v(0, \xi ' )=0 $, $v (y_d , \xi ')\to 0$ as $y_d \to \infty $.
Here we assume that $p (y_d , \xi ') $ decays exponentially as $y_d \to \infty $ and is homogeneous of the generalized degree $s$.
Extending $v$ and $p$ to be zero in $-\infty < y_d <0$, we have 
\begin{equation*}
v(y_d , \xi ' )= \frac{1}{2 |\xi '|_{\Gamma}}   \left( \int_{0}^{y_d } e^{-|\xi '|_{\Gamma}  (y_d - \eta )} p(\eta ,\xi ') d \eta + \int_{y_d }^{\infty} e^{-|\xi '|_{\Gamma}  ( \eta -y_d )} p( \eta ,\xi ') d \eta \right) . 
\end{equation*}
Then, putting $\tau = t\eta $, we have
\begin{gather*}
\begin{split}
&v(t^{-1} y_d , t \xi ' ) \\
&= \frac{t^{s-2} }{2 |\xi '|_{\Gamma}}  \left( \int_{0}^{ y_d } e^{-|\xi '|_{\Gamma}  (y_d - \tau )} p( \tau ,\xi ') d \tau + \int_{ y_d }^{\infty} e^{-|\xi '|_{\Gamma}  ( \tau -y_d )} p( \tau ,\xi ') d \tau \right) \\
&= t^{s-2} v( y_d , \xi ') ,
\end{split}
\end{gather*}
which shows that $v$ is homogeneous of the generalized degree $s-2 $ with respect to $y_d $ and $\xi ' $.
In view of Proposition \ref{S2_prop_degree}, we have $\widetilde{A}_{k,1} E_{k,0} $ is homogeneous of the generalized degree $1$.
Therefore, we obtain $E_{k,1} $ is homogeneous of the generalized degree $-1 $.
Repeating the similar argument inductively, we can show that $E_{k,m} $ is homogeneous of the generalized degree $-m$.
\qed

\medskip

Let $ \beta (\xi ' ) \in C^{\infty} ({\bf R}^{d-1} )$ vanish in a neighborhood of $ \xi ' =0 $, and be equal to one outside a large neighborhood of $\xi ' =0 $. 
Taking $ \psi \in H^{3/2} ( \partial \Omega^0 ) $ with a compact support in $\partial \Omega ^0 $, we define for $y' \in \partial \Omega ^0 $ 
\begin{gather}
\begin{split}
&(Q_{k,m} \psi ) (z'; y' , y_d ) \\
&= (2\pi )^{-(d-1)} \int e^{i y' \cdot \xi '} \beta (\xi ') E_{k,m} (z ' ; y_d , \xi ') \int e^{-i w' \cdot \xi '} \psi ( w' ) dw' d\xi ' ,
\end{split}
\label{S2_def_Gkm}
\end{gather}
and we put 
\begin{equation}
R_{k,N} = \sum_{m=0}^N Q_{k,m} .
\label{S2_def_RkN}
\end{equation}
Letting 
\begin{equation}
 q_{k,m} ( z' ; y' ,y_d ) = (2\pi )^{-(d-1)  } \int e^{i y' \cdot \xi '} \beta (\xi ') E_{k,m} (z ' ; y_d , \xi ') d \xi ' ,
\label{S2_def_Ekm}
\end{equation}
we have that $q_{k,m} $ is a distribution in $\mathcal{S}' $, and 
\begin{gather}
( Q_{k,m } \psi )(z' ; y' ,y_d )= \int q_{k,m} ( z' ; y' -w' , y_d ) \psi ( w') dw ' , \label{S2_eq_Gkmconv} \\
(R_{k,N} \psi )(z' ; y' , y_d )= \int r_{k,N} (z' ; y' -w' , y_d ) \psi (w') dw' , \label{S2_eq_RkNconv}
\end{gather}
with 
$$
 r_{k,N} (z' ; y' - w' , y_d ) = \sum_{m=0}^N  q_{k,m} (z' ; y' - w' , y_d ) .
$$
We represent $A_k $ in the form 
\begin{gather*}
\begin{split}
A_k = &A_{k,0}(z';i\partial_{y'},i\partial_d)+A_{k,1}(z';y'-z',y_d,i\partial_{y'},i\partial_d)\\
&+\sum_{m=2}^NA_{k,m}(\lambda,z';y'-z',y_d,i\partial_{y'},i\partial_d)+A'_{k,N+1}(\lambda,z';y'-z',y_d,i\partial _{y'},i\partial_d) .
\end{split}
\end{gather*}
In the following, we consider  
\begin{gather}
\begin{split}
& A_k r_{k,N} \\
=& \sum_{J = 0}^{N}\sum_{l+m = J}A_{k,l}q_{k,m}+\sum_{J = N+1}^{2N}\sum_{l,m \leq N , l+m = J}A_{k,l}q_{k,m}+A'_{k,N+1}r_{k,N}.
\end{split}
\label{S2_eq_AR}
\end{gather}

\begin{lemma}
Let $l $, $m $ and $N$ be sufficiently large.
We have $ A_{k,l} q_{k,m} \in H^{\gamma }  (\Omega ) $ and $ A'_{k,N+1} r_{k,N} \in H^{ \gamma ' }  (\Omega )$ where $ \gamma = O(l+m)$ and $\gamma ' = O(N) $. 

\label{S2_lem_regurality}
\end{lemma}

Proof.
Note that $A_{k,l}(\lambda,z';y'-z',y_d,i\partial_{y'},i\partial_d)$ and $A'_{k,N+1}(\lambda,z';y'-z',y_d,i\partial_{y'},i\partial_d)$ are operators which are given by sums of terms like $(y' -z')^{\alpha '}  y_d^{\alpha_d } \partial_{y'}^{\beta' } \partial_d^{\beta _d} $ up to a smooth function with $-| \alpha '| - \alpha_d + |\beta '| + \beta_d =2-l $ or $2-(N+1)$ and $ |\beta '| + \beta_d \leq 2 $. 
In view of Proposition \ref{S2_prop_degree}, it is sufficient to show 
\begin{equation}
 (y' ) ^{ \alpha ' } y_d^{\alpha_d} q_{k,m} (z; y'  , y_d ) \in H^{\gamma } (\Omega ),
\label{S2_eq_Hlm}
\end{equation}
since the derivative $\partial_{y'}^{\beta'}\partial_d^{\beta_d}$ is order zero, one or two.

Now we have 
\begin{gather*}
\begin{split}
& (y' ) ^{ \alpha ' } y_d^{\alpha_d} q_{k,m} (z; y'  , y_d )  \\
&= i^{|\alpha'|} (2\pi )^{-(d-1)}  \int e^{i y' \cdot \xi '} \partial _{\xi '}^{ \alpha '  } \big( y_d^{\alpha_d} \beta (\xi ') | \xi '|^{-m}  E_{k,m} (z' ; |\xi '|y_d , \xi  '  / |\xi '| ) \big) d\xi ' .
\end{split}
\end{gather*}
Since $y_d^{\alpha_d} | \xi '|^{-m}  E_{k,m} (z' ; |\xi '|y_d , \xi  '  / |\xi '| )  $ is homogeneous of the generalized degree $-m-\alpha_d $, using proposition \ref{S2_prop_degree}, we have 
$$
\left|\partial_{\xi '}^{ \alpha'}\big(y_d^{\alpha_d}\beta(\xi')E_{k,m}(z';y_d,\xi')\big)\right| \leq C_{m,\alpha}(1+|\xi'|)^{-m-|\alpha'|-\alpha_d},
$$
which implies (\ref{S2_eq_Hlm}). 
\qed

\begin{theorem}
Let $N>1 $ be sufficiently large.
The operator $ R_{k,N} $ satisfies 
\begin{equation}
A_k R_{k,N} \psi \in H^s (\Omega ) , \quad R_{k,N} \psi \big| _{y_d =0} - \psi \in C^{\infty} (\partial \Omega ^0 ) ,
\label{S2_eq_localreg_thm}
\end{equation}
for $\psi \in H^{3/2} ( \partial \Omega ^0 ) $ which has a compact support in $\partial \Omega ^0$ and $s = O(N)$.

\label{S2_thm_localreg}
\end{theorem}

Proof.
Note that 
\begin{gather}
\begin{split}
&A_{k,l}(\lambda,z;y'-z',y_d,i\partial_{y'},i\partial_d)q_{k,m}(z';y'-w',y_d)\\
&= (2\pi )^{-(d-1)} \int e^{i (y' -w') \cdot \xi '} \widetilde{A}_{k,l} \big( \beta (\xi ') E_{k,m} ( z' ; y_d , \xi ') \big) d\xi ' .
\end{split}
\label{S2_eq_Fourier1}
\end{gather}
Summing up both sides of (\ref{S2_eq_sys0})-(\ref{S2_eq_sys2}), we have 
\begin{equation}
\sum_{J=0}^N \sum_{l+m=J} \widetilde{A}_{k,l} E_{k,m} (z' ; y_d , \xi ') =0 .
\label{S2_eq_sys33}
\end{equation} 
In view of Lemma \ref{S2_lem_regurality} and (\ref{S2_eq_AR}), we have that (\ref{S2_eq_Fourier1}) and (\ref{S2_eq_sys33}) imply $A_k R_{k,N} \psi \in H^s (\Omega )$ for $s =O(N) $.

We have that
\begin{gather*}
\begin{split}
& R_{k,N} \psi (y' , y_d) - \psi (y') \\
&= (2\pi )^{-(d-1)} \iint e^{i (y' -w') \cdot \xi '} \left( \sum_{m=0}^N \beta (\xi ')  E_{k,m} (z' ; y_d , \xi ') -1 \right) \psi (w') d\xi ' dw' \\
& \to (2\pi )^{-(d-1)} \iint e^{i (y' -w') \cdot \xi '} \left(  \beta (\xi ')   -1 \right) \psi (w') d\xi ' dw'  ,
\end{split}
\end{gather*}
as $y_d \to 0 $.
Since $ \beta (\xi ') -1 \in C_0^{\infty} ({\bf R}^{d-1} )$, we have $ R_{k,N} \psi \big| _{y_d =0} - \psi (y') \in C^{\infty} (\partial \Omega ^0 )$.
\qed

\medskip

\textit{Remark.} 
The formal sum 
$$
(R_k \psi )( z' ; y' , y_d )  = \int \sum_{m=0}^{\infty} q_{k,m} (z' ; y'-w' , y_d ) \psi (w') dw' , 
$$
is a pseudo-differential operator (see \cite{VaGu}).
In general, a linear operator $P$ on a $d$-dimensional compact manifold $M$ is a pseudo-differential operator of order $l$ if there exist homogeneous functions $ p_j (x,\xi ) \in C^{\infty} (M, {\bf R}^d / \{ 0\} )$ in $\xi $ with homogeneous degree $l-j$ such that for a function $u$ with support in a local coordinate neighborhood $U \subset M $,
$$
Pu (x)= (2\pi )^{-d} \iint e^{i (x -y)\cdot \xi } \beta (\xi ) \sum_{j=0}^N p_j (x,\xi ) u(y) dyd\xi + T_{N+1} u , \quad x\in U,
$$
where $\beta \in C^{\infty} ({\bf R}^d )$ is an arbitrary function which satisfies $ \beta (\xi )=0$ for $|\xi |\leq 1 $ and $\beta (\xi ) =1 $ for $|\xi | \geq 2 $, and $T_{N+1} $ is an operator which increases the smoothness i.e. $H^s (M) \to H^{s+O(N) } (M) $ for any $s\in {\bf R} $. 
The principal symbol of $P$ is $p_0 (x,\xi )$ and the full symbol of $P$ is the formal sum $\sum _j p_j (x,\xi )$.
Then the ellipticity of $P$ is defined by $p_0 (x,\xi ) \not= 0$ for all $ \xi \not= 0$.
Here this means that we can construct the parametrix of $P$ (see \cite{Ho}).
Therefore, if $P$ is an elliptic pseudo-differential operator, $P$ is Fredholm.

\medskip

Since we have $ \partial _{\nu_k} = - \partial _d $ in $y$-coordinates, we can show the following fact. 
As a consequence of Corollary \ref{S2_cor_singular_sol2} and Theorem \ref{S2_thm_localreg}.
See also Lemma 11 and Theorem 14 in \cite{VaGu}.

\begin{cor}

(1) When $\lambda $ is not a pole of $\Lambda_k (\lambda )$, $\Lambda_k (\lambda )$ is a pseudo-differential operator on $H^{3/2} ( \Gamma )$ with the full symbol given by the following asymptotic series :
\begin{equation}
\Lambda_k ( \lambda ; y' ,\xi ' )= - \sum_{m=0}^{\infty} \partial_d E_{k,m} ( y' ; y_d ,\xi ') \Big| _{y_d =0 } , \quad y' \in \partial \Omega ^0  .
\label{S2_eq_DNsymbol}
\end{equation}
(2) When $\lambda = \lambda_0 $ is a pole of $\Lambda_k (\lambda  ) $, the regular part $H_k (\lambda )$ of $\Lambda_k (\lambda )$ at $\lambda_0 $ is a pseudo-differential operator on $B_k (\lambda _0 )^c$ with the full symbol given by (\ref{S2_eq_DNsymbol}).
\label{S2_cor_DNsymbol}
\end{cor}

%
%

\subsection{Principal symbol of the D-N map}

We compute the principal symbol of $\Lambda_1 (\lambda )- \Lambda_2 (\lambda )$.
In the following, we denote by $ \partial _{\nu_k}^m$ for $m\geq 1$ higher order normal derivatives on $\Gamma $ associated with $M_k $.
In $y$-coordinates, we can locally represent $\partial _{\nu_k}^m = (-1)^m \partial_d^m $.
Under the assumption (A-1), we additionally assume on $\Gamma$ that

\medskip

{\bf (A-2)} The metrics $g_1 $, $g_2 $ and the indices of refraction $n_1$, $n_2 $ satisfy one of following two cases :

\medskip

{\bf (A-2-1)} For all $x\in \Gamma $, $ \partial _{\nu_1}^m g_1^{ij} (x) = \partial _{\nu_2}^m g_2^{ij} (x) $ for $m\leq 2$, $i,j=1,\cdots ,d$, and  $n_1 (x) \not= n_2 (x)$,

\medskip

or

\medskip

{\bf (A-2-2)}  For all $x\in \Gamma$, $ \partial _{\nu_1}^m g_1^{ij } (x) = \partial _{\nu_2}^m g_2^{ij} (x) $ for $m\leq 3$, $i,j=1,\cdots ,d$, and $n_1 (x) = n_2 (x) $, $ \partial _{\nu_1} n_1 (x) \not= \partial _{\nu_2} n_2 (x) $.

\medskip

Note that, under the assumptions (A-1) with (A-2-1) or (A-2-2), we can see $\widetilde{A}_{1,m} = \widetilde{A}_{2,m} $ for $ m\leq 1 $ or $m\leq 2 $, respectively.

When $\lambda = \lambda_0 $ is a pole of $\Lambda_1 (\lambda ) -\Lambda_2 (\lambda )$, we define a subspace $B (\lambda _0 ) $ of $H^{3/2} (\Gamma )$ by $B( \lambda_0 )= \widetilde{B}_1 (\lambda_0 ) \cup \widetilde{B}_2 (\lambda_0 )$ where 
$\widetilde{B}_k (\lambda_0 ) = B_k (\lambda _0 ) $ if $\lambda_0 $ is a Dirichlet eigenvalue of $-\Delta_{g_k} -\lambda n_k$, and $\widetilde{B}_k (\lambda_0 ) = \emptyset $ if otherwise. 
We denote by $B(\lambda_0 )^c $ the orthogonal complement of $B(\lambda_0 )$ in $L^2 (\Gamma )$.

When $\lambda = \lambda_0 $ is a pole of $\Lambda_1 (\lambda ) - \Lambda_2 (\lambda )$, we call $\Lambda_1 (\lambda ) - \Lambda_2 (\lambda )$ Fredholm if its regular part $H _{\lambda_0} (\lambda )$ is Fredholm.

\begin{lemma} 
In the following, we suppose $\lambda \not= 0$. \\
(1) Let $\lambda $ be not a pole of $\Lambda_1 (\lambda ) - \Lambda_2 (\lambda )$.
For the case (A-2-1), we have $ \Lambda _1 (\lambda ) - \Lambda _2 (\lambda ) : H^{3/2} (\Gamma ) \to H^{5/2} (\Gamma )$ is an elliptic pseudo-differential operator with the principal symbol
\begin{equation}
- \frac{ \lambda    ( n_1 (x ) - n_2 (x ) ) }{ 2| \xi ' | _{\Gamma} } , \quad x\in \Gamma , \quad \xi' \in {\bf R}^{d-1} . 
\label{S2_eq_principalsymb1}
\end{equation}
(2) Let $\lambda $ be not a pole of $\Lambda_1 (\lambda ) - \Lambda_2 (\lambda )$.
For the case (A-2-2), we have $ \Lambda _1 (\lambda ) - \Lambda _2 (\lambda ) : H^{3/2} (\Gamma ) \to H^{7/2} (\Gamma )$ is an elliptic pseudo-differential operator with the principal symbol
\begin{equation}
 \frac{ \lambda  ( \partial_{\nu_1} n_1 (x ) - \partial_{\nu_2} n_2 (x ) )}{ 4| \xi ' | _{\Gamma} ^2 } , \quad x\in \Gamma , \quad \xi' \in {\bf R}^{d-1} . 
\label{S2_eq_principalsymb2}
\end{equation}
(3) When $\lambda $ is a pole of $\Lambda_1 (\lambda ) - \Lambda_2 (\lambda )$, the regular part of $ \Lambda_1 (\lambda ) - \Lambda_2 (\lambda )$ is pseudo-differential operator on $B(\lambda_0 )^c $ with order $-1 $ for (A-2-1) or $-2 $ for (A-2-2).
Its principal symbol is given by (\ref{S2_eq_principalsymb1}) or (\ref{S2_eq_principalsymb2}), respectively. \\ 
(4) For both of (A-2-1) or (A-2-2), $\Lambda_1 (\lambda )- \Lambda_2 (\lambda )$ is Fredholm for $\lambda \in {\bf C} \setminus \{ 0\} $. 

\label{S2_lem_principalsymb}
\end{lemma}

Proof.
Let $n_1 $ and $n_2 $ satisfy (A-2-1).
In $y$-coordinates, we have $ \widetilde{A}_{1,0} = \widetilde{A}_{2,0} $, $\widetilde{A}_{1,1} = \widetilde{A}_{2,1} $ and $\widetilde{A}_{1,2} - \widetilde{A}_{2,2} = -\lambda (n_1 (y',0) - n_2 (y' ,0) ) $.
Then $ E_{1,0} = E_{2,0} = e^{-|\xi '|_{\Gamma} y_d } $, $ E_{1,1} = E_{2,1} $ and 
$$
(-\partial_d^2 +|\xi '|^2_{\Gamma} ) (E_{1,2} -E_{2,2} )= \lambda (n_1 (y' ,0) - n_2 (y' ,0) ) e^{-|\xi '|_{\Gamma} y_d } .
$$
A particular solution of this equation is 
$$
\frac{ \lambda ( n_1 (y',0 )-n_2 (y',0) )}{2 |\xi '|_{\Gamma} } y_d e^{-|\xi '|_{\Gamma} y_d } ,
$$
which vanishes at $y_d =0 $ and $y_d \to \infty $.
Then we can take it as $E_{1,2} - E_{2,2} $, and $-\partial_d (E_{1,2} - E_{2,2} ) $ at $y_d =0 $ is the principal symbol of $\Lambda_1 (\lambda )- \Lambda_2 (\lambda )$.
In view of the assertion (1) in Corollary \ref{S2_cor_DNsymbol}, we have the assertion (1).

Next we assume that $n_1 $ and $n_2 $ satisfy (A-2-2).
As above, we have $ \widetilde{A}_{1,j} = \widetilde{A}_{2,j} $ for $ j=0,1,2$, and $ \widetilde{A}_{1,3} - \widetilde{A}_{2,3} = -\lambda ( \partial_d n_1 (y',0 ) -\partial_d n_2 (y',0 )) y_d $.
Then we have 
$$
E_{1,3} - E_{2,3} = \frac{\lambda}{4}  ( \partial_d n_1 (y',0 ) - \partial_d n_2 (y', 0) ) \frac{y_d}{|\xi '|_{\Gamma} } \left( y_d + \frac{1}{|\xi ' |_{\Gamma} } \right) e^{-| \xi ' |_{\Gamma} y_d } .
$$
Hence we obtain the assertion (2). 

In view of Corollary \ref{S2_cor_singular_sol2} and the assertion (2) in Corollary \ref{S2_cor_DNsymbol}, we can show the assertion (3) by the similar way.

The ellipticity of $\Lambda_1 (\lambda )-\Lambda _2 (\lambda ) $ implies that $ \Lambda_1 (\lambda )-\Lambda _2 (\lambda ) $ is Fredholm for $\lambda \in {\bf C} \setminus \{ 0 \} $.
\qed

%
%
%

\section{Interior transmission eigenvalues}

Let us list our assumptions again : 

\medskip

{\bf (A-1)} $M_1 $ and $M_2 $ have a common boundary $\Gamma := \partial M_1 = \partial M_2 $. 
$\Gamma $ is a disjoint union of a finite number of connected and closed components. 
The metrics satisfy $g_1 =g_2 $ on $\Gamma $.

\medskip

{\bf (A-2)} The metrics $g_1 $, $g_2 $ and the indices of refraction $n_1$, $n_2 $ satisfy one of following two cases :

\medskip

{\bf (A-2-1)} For all $x\in \Gamma $, $ \partial _{\nu_1}^m g_1^{ij} (x) = \partial _{\nu_2}^m g_2^{ij} (x) $ for $m\leq 2$, $i,j=1,\cdots ,d$, and  $n_1 (x) \not= n_2 (x)$,

\medskip

or

\medskip

{\bf (A-2-2)}  For all $x\in \Gamma$, $ \partial _{\nu_1}^m g_1^{ij} (x) = \partial _{\nu_2}^m g_2^{ij } (x) $ for $m\leq 3$, $i,j=1,\cdots ,d$, and $n_1 (x) = n_2 (x) $, $ \partial _{\nu_1} n_1 (x) \not= \partial _{\nu_2} n_2 (x) $.

\medskip

Throughout of \S 3, we suppose the above conditions.

\subsection{Discreteness of the set of ITEs}

For the proof of discreteness, we need to use the analytic Fredholm theory which was generalized by \cite{Bl}.
See also Appendix A in \cite{Sk}.
Let $H_1 $ and $H_2 $ are Hilbert spaces.
We take a connected open domain $D\subset {\bf C} $.
An operator valued function $A(z) : H_1 \to H_2 $ for $z\in D $ is finitely meromorphic if the principal part of the Laurent series at a pole of $A(z)$ is a finite rank operator.
In particular, $\Lambda_k ( \lambda ) : H^{3/2} (\Gamma ) \to H^{1/2} (\Gamma ) $ is finitely meromorphic in ${\bf C}\setminus \{ 0\} $ as has been seen in Proposition \ref{S2_prop_green_meromor}.

\begin{theorem}

Suppose an operator valued function $A(z) : H_1 \to H_2 $, $z\in D$, is finitely meromorphic and Fredholm.
If there exists its bounded inverse $A(z_0 )^{-1} : H_2 \to H_1 $ at a point $z_0 \in D$, then $A( z)^{-1} $ is finitely meromorphic and Fredholm in $D$. 

\label{S3_thm_analyticfredholm}
\end{theorem}

From the above theorem, if $\Lambda_1 (\lambda )-\Lambda_2 (\lambda ) $ is invertible at a point $ \lambda\in {\bf C} \setminus \{ 0\} $, $ \Lambda_1 (\lambda )- \Lambda _2 (\lambda )$ is invertible in ${\bf C} \setminus (\{ 0\} \cup S' )$ for a discrete subset $S'$ of ${\bf C} $. 
Therefore, for the proof of the discreteness, we have only to show that  $ \Lambda_1 (\lambda )- \Lambda _2 (\lambda )$ is invertible for some $\lambda \in {\bf C} \setminus \{ 0 \}$.

We expand the symbol of $A_k$ centered at $(z' ,0)\in \partial \Omega^0 $ by the same manner in \S 2.2. 
However, here we change the definition of homogeneous functions with generalized degree $s$ by
\begin{equation}
F( t \kappa ; t^{-1} y' , t^{-1} y_d , t \xi ' , t \xi _d ) = t^s F( \kappa ;  y' , y_d ,  \xi ' ,  \xi _d ) , \quad t> 0 , \quad \kappa = \sqrt{\lambda},
\label{S3_def_parameterhomogeneous}
\end{equation}
for $\lambda \in {\bf C} \setminus \{ 0 \} $, taking a suitable branch of $ \kappa = \sqrt{\lambda } $.
We gather terms of the same generalized degree in the sense (\ref{S3_def_parameterhomogeneous}), and we denote the symbol in $y$-coordinates by
$$
A_k ( \kappa ;y' , y_d , \xi ' , \xi_d )= \sum_{m=0}^N \mathcal{A}_{k,m} ( \kappa , z' ; y ' -z' , y_d , \xi ' , \xi_d ) , 
$$
up to the remainder term where $\mathcal{A} _{k,m} $ is homogeneous of degree $2-m $.
In particular, putting $ \widetilde{\mathcal{A}}_{k,m}^{(\lambda )} = \mathcal{A}_{k,m} ( \kappa , z' ; -i \partial _{\xi'} , y_d , \xi ' , i\partial_d )$, we have 
\begin{gather}
\widetilde{\mathcal{A}}_{k,0}^{(\lambda )} = -\partial_d^2 + |\xi '| ^2 _{\Gamma} -\lambda n_k (z' , 0) , \label{S3_eq_Ak0lambda} \\
\widetilde{\mathcal{A}}_{k,1}^{(\lambda )} = \widetilde{A}_{k,1} + \lambda  \widetilde{B} _{k,1}^{(\lambda )} , \label{S3_eq_Ak1lambda} 
\end{gather}
where $\widetilde{A}_{k,1} $ is defined by (\ref{S2_def_Atilde2}) and 
$$
\widetilde{B}_{k,1}^{(\lambda )} = i  \nabla _{y'} n_k (z' ,0) \cdot \nabla _{\xi '}  - y_d \partial_d n_k (z' ,0)  . 
$$
We denote by $\{ E_{k,m}^{(\lambda  )} \} _{m\geq 0 } $ the solution of 
\begin{gather}
\widetilde{\mathcal{A}}_{k,0}^{(\lambda )} E_{k,0}^{(\lambda )} = 0 , \label{S3_eq_system1} \\
\widetilde{\mathcal{A}}_{k,0}^{(\lambda )} E_{k,m}^{(\lambda )} = -\sum _{n=0}^m \widetilde{\mathcal{A}}_{k,n}^{(\lambda )} E_{k,m-n} ^{(\lambda )} ,   \quad m\geq 1 ,
\label{S3_eq_system2}
\end{gather}
with the boundary condition $ E_{k,0}^{(\lambda )} \big| _{y_d =0 } =1 $, $E_{k,m}^{(\lambda )} \big| _{y_d =0 } =0 $ for $m\geq 1 $ and $E_{k,m}^{(\lambda )} \to 0 $ as $y_d \to \infty $ for $m \geq 0 $.

In order to apply the theory of parameter-dependent pseudo-differential operators to $\Lambda_1 (\lambda )-\Lambda_2 (\lambda )$, we recall its definition. 
Let $M$ be a $d$-dimensional compact manifold without boundary.
We call $ p(x,\xi , \tau ) \in C^{\infty} (M \times {\bf R}^d \times {\bf R} _{\geq 0} )$ a \textit{uniformly estimated polyhomogeneous symbol of order $s  $ and regularity $r $} if $p(x,\xi ,\tau )$ satisfies 
\begin{gather}
\begin{split}
& | \partial_x^{\alpha} \partial _{\xi}^{\beta} \partial _{\tau}^j p(x,\xi ,\tau ) | \\
& \leq C _{\alpha \beta j} \left( \langle \xi \rangle ^{r-|\beta |} + (|\xi |^2 +\tau^2 +1) ^{(r-|\beta |)/2}  \right)  (|\xi |^2 +\tau^2 +1) ^{(s-r-j)/2} ,
\end{split}
\label{S3_def_phomsymbol1}
\end{gather}  
on $ M \times {\bf R}^d \times {\bf R} _{\geq 0}$ for constants $C_{\alpha \beta j} >0$, and $p(x,\xi ,\tau )$ has the asymptotic expansion
\begin{equation}
p (x,\xi ,\tau ) \sim \sum _{l=0} ^{\infty} p_{s-l} (x,\xi ,\tau ) ,
\label{S3_def_phomsymbol2}
\end{equation} 
where $p_{s-l} (x,\xi ,\tau )$ is homogeneous with generalized degree $s-l $ with respect to $\xi,\tau $ in the sense of 
\begin{equation}
p_{s-l} ( x ,t\xi ,t\tau )= t^{s-l}p_{s-l} (x,\xi ,\tau ) ,\quad t>0 .
\label{S3_def_hompara}
\end{equation} 
A pseudo-differential operator $P(\tau )$ on $M$ with a uniformly estimated polyhomogeneous symbol $p(x,\xi ,\tau )$ is said to be \textit{uniformly parameter elliptic} if the principal symbol $p_{d} (x,\xi ,\tau )  $ does not vanish when $|\xi| +\tau \not= 0$.
For more information and general theory on parameter-dependent operators, one can refer Chapters 2 and 3 in \cite{Gu}.

Let us turn to $\Lambda_1 (\lambda )-\Lambda_2 (\lambda )$.
For $ \lambda \in {\bf C} \setminus {\bf R} _{\geq 0} $, we put $\sqrt{\lambda} = \tau e^{i \theta} $ with $\tau >0$ and $ \theta \in {\bf R} $ such that $ \theta \not= 0$ modulo $\pi $. 
In the following, we fix a suitable $\theta $ and put 
\begin{equation}
R(\tau ) = \tau^{-2} e^{-2i\theta} (\Lambda_1 (\tau^2 e^{2i\theta} ) - \Lambda_2 (\tau^2 e^{2i\theta} )) .
\label{S3_def_Rtau}
\end{equation}

\begin{lemma}
Let $\lambda = \tau^2 e^{2i \theta} \in {\bf C} \setminus {\bf R} _{\geq 0} $. \\
(1) We assume that (A-2-1) holds.
Then $R( \tau ) $ is uniformly parameter elliptic with order $-1$ and regularity $ \infty $.
Its principal symbol is 
\begin{equation}
\frac{ -  (n_1 (x) - n_2 (x))}{\sqrt{ |\xi '|^2 _{\Gamma} - \tau^2 e^{2i\theta} n_1 (x)} + \sqrt{ |\xi '|^2 _{\Gamma} - \tau^2 e^{2i \theta} n_2 (x)}} , \quad x\in \Gamma , \quad \xi ' \in {\bf R}^{d-1} .
\label{S3_eq_symbola31}
\end{equation}
(2) We assume that (A-2-2) holds.
Then $R( \tau ) $ is uniformly parameter elliptic with order $-2$ and regularity $\infty $.
Its principal symbol is 
\begin{equation}
\frac{  (\partial _{\nu_1} n_1 (x) - \partial _{\nu_2} n_2 (x))}{ 4( |\xi '|^2_{\Gamma} - \tau^2 e^{2i\theta} n(x))}  , \quad x\in \Gamma , \quad \xi ' \in {\bf R}^{d-1} ,
\label{S3_eq_symbola32}
\end{equation}
where $n(x ) := n_1 (x) = n_2 (x)$.

\label{S3_lem_parametersymbol}
\end{lemma}

Proof. 
We fix an arbitrary point $ (z' ,0) \in \partial \Omega ^0 $.
Suppose that (A-2-1) holds.
Obviously we have 
$$
E_{k,0}^{(\lambda )} ( z' ; \xi ' ,y_d ) = \mathrm{exp} \left( -\sqrt{ | \xi ' | _{\Gamma}^2 - \lambda n_k (z',0 )} y_d \right) .
$$
Under the assumption, we also have $ \widetilde{\mathcal{A}}_{1,0}^{(\lambda )} \not= \widetilde{\mathcal{A}} _{2,0} ^{(\lambda )} $ so that $ E_{1,0}^{(\lambda )} \not= E_{2,0} ^{(\lambda )} $.
Then the principal symbol $- \partial_d ( E_{1,0}^{(\lambda )} - E_{2,0} ^{(\lambda )}) \big| _{y_d =0 }$ of $\Lambda_1 (\lambda )  -\Lambda_2 (\lambda ) $ is given by 
$$
\frac{ - \lambda  (n_1 (x) - n_2 (x))}{\sqrt{ |\xi '|^2 _{\Gamma} - \lambda n_1 (x)} + \sqrt{ |\xi '|^2 _{\Gamma} - \lambda n_2 (x)}} .
$$
This shows (\ref{S3_eq_symbola31}).

Let us consider the case (A-2-2).
In view of $n_1 = n_2 (=n) $ on $\Gamma$, we have $\widetilde{\mathcal{A}}_{1,0}^{(\lambda )} = \widetilde{\mathcal{A}} _{2,0} ^{(\lambda )} $ so that 
$$
 E_{1,0}^{(\lambda )} (z' ; \xi ' , y_d ) = E_{2,0} ^{(\lambda )} (z' ; \xi ' ,y_d ) = \mathrm{exp} \left( - \sqrt{ | \xi ' | _{\Gamma}^2 - \lambda n (z',0 )} y_d \right) .
$$ 
Since we have assumed (A-1) and (A-2-2), we have
$$
\widetilde{\mathcal{A}} _{1,1}^{(\lambda )} - \widetilde{\mathcal{A}} _{2,1}^{(\lambda )} = -\lambda ( \partial_d n_1 (z',0 )- \partial_d n_2 (z',0 )) y_d .
$$
Then $E_{1,1}^{(\lambda )} - E_{2,1}^{(\lambda )} $ satisfies the equation 
\begin{gather*}
\begin{split}
& (-\partial_d^2 + |\xi '| _{\Gamma}^2 -\lambda n(z',0 )) (E_{1,1}^{(\lambda )} - E_{2,1}^{(\lambda )} ) \\
&= \lambda ( \partial_d n_1 (z',0 )- \partial_d n_2 (z',0 )) y_d \,  \mathrm{exp} \left(- \sqrt{ | \xi ' | _{\Gamma}^2 - \lambda n (z',0 )} y_d \right) .
\end{split}
\end{gather*}
Precisely, we obtain 
\begin{gather*}
\begin{split}
&E_{1,1}^{(\lambda )} (z' ; \xi ' , y_d )- E_{2,1}^{(\lambda )} (z' ; \xi ' , y_d ) = -\frac{\lambda}{4} ( \partial_d n_1 (z',0 )- \partial_d n_2 (z',0 )) \\
& \cdot \left( \frac{y_d^2}{ \sqrt{ | \xi ' | _{\Gamma}^2 - \lambda n (z',0 )}} + \frac{y_d}{ | \xi ' | _{\Gamma}^2 - \lambda n (z',0 )} \right) \mathrm{exp} \left(- \sqrt{ | \xi ' | _{\Gamma}^2 - \lambda n (z',0 )} y_d \right) .
\end{split}
\end{gather*}
Then the principal symbol $- \partial_d ( E_{1,1}^{(\lambda )} - E_{2,1} ^{(\lambda )}) \big| _{y_d =0 }$ of $\Lambda_1 (\lambda )  -\Lambda_2 (\lambda ) $ is given by
$$
\frac{ \lambda  (\partial _{\nu_1} n_1 (x) - \partial _{\nu_2} n_2 (x))}{ 4( |\xi '|^2_{\Gamma} - \lambda  n(x))} .
$$ 
This shows (\ref{S3_eq_symbola32}).
\qed

\medskip

In view of Lemma \ref{S3_lem_parametersymbol}, we can obtain a uniform estimate in $\tau $ of $ R(\tau ) $ and its inverse.
In the following, we define the Hilbert space $H^{m,t} (\Gamma )$ for $t \ge 1$ by the norm 
$$
\| f\|^2_{H^{m,t} (\Gamma )} = \| f\|^2 _{H^m (\Gamma )} + t^{2m} \| f\|^2 _{L^2 (\Gamma )} .
$$

\begin{lemma}
For sufficiently large $\tau >0$, there exists $ R(\tau ) ^{-1} : H^{m, \tau } (\Gamma ) \to H^{m-s , \tau } (\Gamma )$ for any $m\in {\bf R}$ where $s=1$ for (A-2-1) or $s=2$ for (A-2-2).
\label{S3_thm_DNinverse}
\end{lemma}

Proof.
In view of Lemma \ref{S3_lem_parametersymbol}, we can construct the parametrix of $R(\tau )$. 
The theorem is a direct consequence of Theorem 3.2.11 in \cite{Gu}.
\qed

\medskip

Let us turn to the case $\zeta \not= 0$.
In view of $ \Lambda_1 (\lambda ) -\Lambda_2 (\lambda ) -\zeta = \zeta^{1/2} ( \zeta ^{-1/2} (\Lambda_1 (\lambda ) -\Lambda_2 (\lambda ) ) \zeta ^{-1/2} -1) \zeta ^{1/2}$, we put 
\begin{equation}
K(\lambda )= \zeta ^{-1/2} ( \Lambda_1 (\lambda ) -\Lambda_2 (\lambda ) ) \zeta ^{-1/2}.
\label{S3_def_Klambda}
\end{equation}
Since $\zeta \in C^{\infty} (\Gamma )$ is strictly positive or strictly negative and $\Lambda_1 (\lambda ) -\Lambda_2 (\lambda ) $ has a negative order, the operator $K(\lambda )$ is compact in $L^2 (\Gamma )$ when $\lambda$ is not a pole.
Since $K(\lambda )$ is meromorphic with respect to $\lambda $, we have the following lemma.
The proof is completely same of and 2.4 in \cite{LaVa}.
Note that we will refer the above lemma again later.

\begin{lemma}
Let $ \{ \kappa_j (\lambda ) \} $ be the set of eigenvalues of $K(\lambda )$.
Then every $\kappa_j (\lambda )$ is meromorphic with respect to $\lambda $.
If $\lambda_0$ is a pole of $K(\lambda )$ and $p$ is the rank of the residue of $K(\lambda )$ at $\lambda_0 $, $p$ eigenvalues and its eigenfunctions have a pole at $\lambda_0 $.
Moreover, $\mathrm{res} _{\lambda  = \lambda_0 } \kappa_j (\lambda ) $ are eigenvalues of $\mathrm{res} _{\lambda = \lambda_0 } K(\lambda )$.
\label{S3_lem_analyticev}
\end{lemma}

As a consequence, we have the following lemma.

\begin{lemma}
There exist $\lambda \in {\bf C} \setminus {\bf R} _{\geq 0} $ such that $1\not\in \{ \kappa_j (\lambda ) \} $. 
In particular, $K(\lambda )-1$ has the bounded inverse for some $ \lambda \in {\bf C}\setminus \{ 0 \} $.

\label{S3_lem_evone}
\end{lemma}

Proof. 
Note that the set $\mathcal{A} = \{ \lambda \in {\bf C} \setminus \{ 0 \} \ ; \ \lambda \text{ is not a pole of } K(\lambda ) \} $ is a connected domain in ${\bf C} \setminus \{ 0 \} $. 
Since $K(\lambda  )$ is compact, $\{ \kappa_j (\lambda ) \} $ is the set of eigenvalues of finite multiplicities with the only possible accumulation point at $0$.

We take a point $ \lambda_1 \in {\bf C} \setminus {\bf R} _{\geq 0}$ such that $ \kappa _j (\lambda _1 )= \cdots = \kappa _{j+l-1} (\lambda_1 )=1 $.
In view of the discreteness of eigenvalues, there exists a small constant $\epsilon_0 >0 $ such that $|  \kappa _m (\lambda_1 ) -1 | > \epsilon_0 $ for $m\not\in \{ j , j+1 , \cdots , j+l-1 \} $.
Taking a sufficiently small $\delta >0 $, we also have $|  \kappa _m (\lambda ) -1 | > \epsilon_0 $ for $| \lambda - \lambda_1 | < \delta $.

Suppose that there exists an eigenvalue $\kappa_{j'} (\lambda )$ with $j'\in \{ j , j+1 , \cdots , j+l-1 \}  $ such that $ \kappa _{j'} (\lambda )=1 $ in $\{ \lambda \in {\bf C} \ ; \  | \lambda - \lambda_1 | <\delta \} $.
Since $ \kappa _{j'} (\lambda )$ is analytic in $ \mathcal{A} $, we have $\kappa_{j'} (\lambda )=1 $ in $\mathcal{A} $. 
We take a pole $ \lambda_0 $ of $ \kappa_{j'} (\lambda )$.
In a small neighborhood of $\lambda_0$, $\kappa _{j'} (\lambda ) $ can be written by 
$$
\kappa _{j'} (\lambda )= \frac{ \mathrm{res} _{\lambda = \lambda_0} \kappa _{j'} (\lambda )}{ \lambda_0 -\lambda } + \widetilde{\kappa} _{j'} (\lambda ) ,
$$
where $ \widetilde{\kappa} _{j'} (\lambda ) $ is analytic in this neighborhood. 
However, we obtain
$$
\mathrm{res} _{\lambda = \lambda_0} \kappa _{j'} (\lambda )= ( \lambda_0 -\lambda )(1-\widetilde{\kappa} _{j'} (\lambda ) ) \to 0 ,
$$
as $\lambda \to \lambda_0 $.
This is a contradiction.
\qed

\medskip

Now we have our first main theorem as a corollary of Theorem \ref{S3_thm_analyticfredholm}, Lemma \ref{S3_thm_DNinverse} and Lemma \ref{S3_lem_evone}.
We take an arbitrary closed sector ${\bf S}_0 $ centered at the origin such that ${\bf S}_0 \cap {\bf R} _{>0} = \emptyset $. 
We put ${\bf S}_0^e : = {\bf S}_0 \cap \{ \lambda \in {\bf C} \ ; \ |\lambda | \geq 1 \} $.

\begin{theorem}
We assume (A-1) and one of (A-2-1) and (A-2-2).
The set of ITEs consists of a discrete subset of ${\bf C}$ with the only possible accumulation points at $0$ and infinity. 
There exist at most finitely many ITEs in ${\bf S}_0^e $.
\label{S3_thm_discreteness}
\end{theorem}

Proof.
Note that $  \Lambda_1 (\lambda ) -\Lambda_2 (\lambda ) - \zeta  $ is finitely meromorphic and Fredholm for $ \lambda \in {\bf C} \setminus \{ 0 \} $. 
Lemma \ref{S3_thm_DNinverse} implies that the bounded inverse $ (\Lambda_1 (\lambda ) -\Lambda_2 (\lambda )  ) ^{-1} $ exists for $\lambda \in {\bf S}_0^e $ with sufficiently large $|\lambda|$. 
Lemma \ref{S3_lem_evone} implies that the bounded inverse $(\Lambda_1 (\lambda ) -\Lambda_2 (\lambda ) -\zeta  ) ^{-1} $ exists for some $ \lambda \in {\bf C} \setminus {\bf R} _{\geq 0} $.
In view of Theorem \ref{S3_thm_analyticfredholm}, we obtain the theorem for both of the cases $\zeta =0$ and $\zeta \not= 0$.
\qed

\subsection{Weyl type estimate for interior transmission eigenvalues}

In the following, we use Weyl's law at infinity for Dirichlet eigenvalues of $- n_k ^{-1} \Delta_{g_k }$ on $M_k $. 
The following fact is a direct consequence of Theorem 1.2.1 in \cite{SaVa}.

\begin{theorem}
Let $ \mathcal{O} _k ( x) = \{ \xi \in {\bf R}^d \ ; \ \sum_{i,j} g^{ij}_k (x) \xi _i \xi _j \leq n_k (x) \} $ for each $x \in M_k$ and 
$$
 v( \mathcal{O}_k  (x) ) := \int _{ \mathcal{O}_k (x) } d\xi ,
$$
be the volume of $\mathcal{O}_k(x)$.
Then $ N_k (\lambda ) := \# \{ j \ ; \ \lambda_{k,j} \leq \lambda \} $ satisfies as $\lambda \to \infty  $ 
\begin{equation}
N_k (\lambda )= V_k \lambda ^{d/2} +O(\lambda ^{(d-1)/2} ) , \quad V_k = (2\pi )^{-d} \int_{M_k} v( \mathcal{O}_k (x) ) dV_k .
\label{S3_eq_weylDirichlet}
\end{equation}

\label{S3_thm_Weyl}
\end{theorem}

Taking an arbitrary point $x^{(0)} \in \Gamma  $, we take a small neighborhood $V\subset \Gamma  $ of $x^{(0)} $ and a sufficiently small open domain $\Omega $ which is diffeomorphic to $ U_1 \cong U_2 $ such that $ \overline{U_1} \cap \Gamma = \overline{U_2} \cap \Gamma =V$ as has been defend in the beginning of \S 2.2. 
Then, identifying $x\in V$ with the corresponding point $y\in \partial \Omega^0 $, we define
\begin{equation}
 \gamma_{0} (x) :=
 \left\{
  \begin{split}
   &\mathrm{sign} ( n_2 ( y) - n_1 (y) ) \quad \text{for} \quad \text{(A-2-1)}, \\
   &\mathrm{sign} (\partial_{\nu_1} n_1(y)-\partial_{\nu_2}n_2(y)) \quad \text{for} \quad \text{(A-2-2)} ,
  \end{split}
 \right.  
\label{S3_eq_sgnn1n2}
\end{equation}
 and
\begin{equation}
\gamma _{ \zeta } (x) :=- \mathrm{sign}(\zeta(y)) \quad \text{for} \quad \zeta \neq 0 .
 \label{S3_eq_sgnzeta}
\end{equation}
Note that $\gamma_0 (x) $ and $ \gamma _{\zeta} (x)$ are well-defined constant functions $\gamma_0 (x) = 1 $ or $-1$ and $\gamma _{\zeta} (x) = 1$ or $-1$ on each connected component of $\Gamma $, respectively.
We also define the function $\gamma$ on $\Gamma$ by
\begin{equation}
 \gamma = 
 \left\{
  \begin{split}
   &\gamma_0 \quad \text{for} \quad \zeta = 0 , \\
   &\gamma_{\zeta } \quad \text{for} \quad \zeta \neq 0.
  \end{split}
 \right.
\label{S3_def_gammadef}
\end{equation}

Generally, the function $\gamma $ can change its value for each connected component.
However, let us impose the following third assumption for the proof of Theorem \ref{S3_thm_WeylNT}.
In the following, we suppose (A-3) for all lemmas.

\medskip

{\bf (A-3)} If $\zeta =0$, then $n_1 (x) -n_2 (x) $ or $ \partial _{\nu_1} n_1 (x) -\partial _{\nu_2} n_2 (x)$ do not change its sign on whole of $\Gamma $. 
If $\zeta \not= 0$, then $\zeta $ does not change its sign on whole of $\Gamma $.
In particular, the function $\gamma  $ is constant $1$ or $-1 $ on $\Gamma$.

\medskip

In the following, we use an auxiliary operator defined by 
\begin{equation}
B(\lambda )= \gamma D_{\Gamma}^{(1 + s)/4} ( \Lambda_1 (\lambda ) -\Lambda_2 (\lambda ) -\zeta) D_{\Gamma}^{(1 + s)/4}.
\label{S3_def_Blambda}
\end{equation}
Here $D_{\Gamma}$ is given by $D_{\Gamma} = -\Delta_{\Gamma}+1$ where $\Delta_{\Gamma} $ is the Laplace-Beltrami operator on $\Gamma $. 
If $\zeta =0 $, we take $s=1$ for (A-2-1) or $s=2$ for (A-2-2).
If $\zeta \not= 0$, we take $s=0$.
Then $B( \lambda )$ is a first order pseudo-differential operator when $\lambda $ is not a pole of $ \Lambda_1 (\lambda ) -\Lambda_2 (\lambda )$.

\begin{lemma}
(1) Suppose $\lambda \not\in \{ \lambda_{1,j} \} \cap \{ \lambda_{2,j} \} $.
Then $\lambda \in {\bf C} $ is an ITE if and only if $\mathrm{Ker} B( \lambda ) \not= \{ 0 \}  $.
The multiplicity of $\lambda $ coincides with $\mathrm{dim}  \mathrm{Ker} B(\lambda ) $. \\
(2) Suppose $\lambda \in \{ \lambda_{1,j} \} \cap \{ \lambda_{2,j} \} $.
Then $\lambda \in {\bf C} $ is an ITE if and only if $ \mathrm{Ker} B(\lambda ) \not= \{ 0 \}  $ or the ranges of $ \gamma D_{\Gamma}^{(1+s)/4} Q_{1,\mathcal{L}(\lambda)} D_{\Gamma}^{(1+s)/4} $ and $ \gamma D_{\Gamma}^{(1+s)/4} Q_{2,\mathcal{L}(\lambda)} D_{\Gamma}^{(1+s)/4} $ have a non trivial intersection. 
The multiplicity of $\lambda $ coincides with the sum of $\mathrm{dim}  \mathrm{Ker} B(\lambda ) $ and the dimension of the above intersection.  
\label{S2_lem_ITEeqconditionB}

\end{lemma}

Proof. 
Since $ -\Delta _{\Gamma} +1 $ is invertible, the lemma is a direct consequence of Lemma \ref{S2_lem_ITEeqcondition}. 
\qed

\medskip

\begin{lemma}
 Let $ \lambda \in {\bf C} \setminus \{ 0\} $ be not a pole of $B(\lambda )$. \\
(1) For $\zeta = 0$, $B(\lambda )$ is a first order, symmetric and elliptic pseudo-differential operator.
Its principal symbol is 
\begin{equation}
 \frac{\lambda\gamma(n_2 (x) -n_1 ( x))}{2}|\xi'|_{\Gamma}, \quad x\in \Gamma , \quad \xi ' \in {\bf R}^{d-1} ,
\label{S3_eq_psymbolB1}
\end{equation}
for (A-2-1), or 
\begin{equation}
 \frac{\lambda\gamma(\partial_{\nu_1}n_1(x)-\partial_{\nu_2}n_2(x))}{4}|\xi'|_{\Gamma}, \quad x\in \Gamma , \quad \xi ' \in {\bf R}^{d-1} ,
\label{S3_eq_psymbolB2}
\end{equation}
for (A-2-2).
 \\
(2) For $\zeta \neq 0$, $B(\lambda)$ is a first order, symmetric and elliptic pseudo-differential operator.
Its principal symbol is 
\begin{equation}
 -\gamma\zeta(x)|\xi'|_{\Gamma}, \quad x \in \Gamma, \quad \xi' \in {\bf R}^{d-1}.
\end{equation}
(3) For $\lambda \in {\bf R} _{>0} $, the spectrum of $B(\lambda )$ is discrete and consists of the set of real eigenvalues $\{ \mu_j (\lambda ) \} $. 
\label{S3_lem_symbB}
\end{lemma}

Proof.
We have the first assertion by direct computation using Lemma \ref{S2_lem_principalsymb}.
From the first assertion, we also see the second assertion.
\qed

\medskip

Since $B(\lambda )$ has a positive principal symbol and $B(\lambda )$ is meromorphic with respect to $\lambda $, we also have the following lemma.
For the proof, see Lemmas 2.3 and 2.4 in \cite{LaVa}.
Note that, in view of (\ref{S2_eq_resDN}), we define the residue $ \mathrm{res} _{\lambda = \lambda_0 } \mu_j (\lambda ) $ of $\mu_j (\lambda )$ at a pole $\lambda_0 $ by the expansion
\begin{equation}
\mu_j (\lambda )= \frac{ \mathrm{res} _{\lambda =\lambda_0 } \mu_j (\lambda )}{\lambda_0 -\lambda } + \widetilde{\mu }_j (\lambda ) ,
\label{S3_def_residue}
\end{equation}
where $\widetilde{\mu}_j (\lambda )$ is analytic in a small neighborhood of $\lambda_0 $.

\begin{lemma}
(1) For each compact interval $I \subset {\bf R}_{ >0} $ such that there is no pole of $B(\lambda )$ in $I$, there exists a constant $C(I) >0 $ such that $\mu_j (\lambda ) \geq -C(I) $ for $\lambda \in I$, $j=1,2, \cdots $.  \\ 
(2) If $B(\lambda )$ is analytic in a neighborhood of $\lambda_0 \in {\bf R} _{>0} $, all eigenvalues $\mu_j (\lambda )$ are analytic in this neighborhood.
If $\lambda_0 \in {\bf R}_{>0} $ is a pole of $B(\lambda )$ and $ p$ is the rank of the residue of $B(\lambda )$ at $\lambda_0 $, $p$ eigenvalues $\mu_j (\lambda )$ and its eigenfunctions have a pole at $\lambda_0 $. 
Moreover, $  \mathrm{res}_{\lambda = \lambda_0 } \mu_j (\lambda ) $ are eigenvalues of $ \mathrm{res}_{\lambda =\lambda_0} B(\lambda )$.

\label{S3_lem_positiveB}
\end{lemma}

We choose a small constant $ \alpha \in (0, \mathrm{min} \{ \lambda _{1,1} , \lambda _{2,1} \} )$.
We define the counting function with multiplicities taken into account : 
\begin{equation}
N_T (\lambda ) = \# \{ j \ ; \ \alpha <\lambda^T_j \leq \lambda \} , 
\label{S3_def_countfunc}
\end{equation}
where $ \lambda_1^T \leq \lambda _2^T  \leq \cdots $ are ITEs included in $(\alpha , \infty )$.

Now we consider the relation between $\{ \lambda_j^T \} $ and $\{ \mu_j (\lambda ) \} $ for $\lambda \in (\alpha , \infty )$.
Roughly speaking, we can evaluate $N_T (\lambda )$ by the number of the singular ITEs and the number of $\lambda$ satisfying $\mu_j(\lambda) = 0$. 
We put 
\begin{equation}
 N_- (\lambda )= \# \{ j \ ; \ \mu_j (\lambda ) <0 \} , \quad \lambda \not\in \{ \lambda_j^T \} \cup \{ \lambda_{1,j} \} \cup \{ \lambda _{2,j} \} .
\label{S3_def_countmu0}
\end{equation}
Assume that $\lambda '$ moves from $\alpha $ to $\infty $.
Since $\mu_j (\lambda ' )$ is meromorphic with respect to $\lambda '$, $N_- (\lambda ')$ changes only when some $ \mu _j (\lambda ')$ pass through $0$ or $\lambda '$ passes through a pole of $B(\lambda ')$.  
When $\lambda '$ moves from $\alpha $ to $\lambda > \alpha $, we denote by $\mathcal{N}_0 (\lambda )$ the change of $N_- (\lambda ) -N_- (\alpha )$ due to the first case, and $\mathcal{N}_{-\infty} (\lambda )$ as the change due to the second case i.e. 
\begin{equation}
 N_- (\lambda )-N_- (\alpha ) = \mathcal{N}_0 (\lambda ) +  \mathcal{N}_{-\infty} (\lambda ) .
\label{S3_eq_diffNB}
\end{equation}
For a pole $ \lambda_0 $ of $B(\lambda )$, we put 
\begin{equation}
\delta \mathcal{N}_{-\infty}(\lambda_0) = N_-(\lambda_0+\epsilon)-N_{-}(\lambda_0-\epsilon),
\label{S3_eq_jumpNinfty}
\end{equation}
for sufficiently small $\epsilon >0 $.

\begin{lemma}
Let $\lambda _0 \in {\bf R}_{>0} $ be a pole of $B(\lambda )$.
We have $ \delta \mathcal{N}_{-\infty} (\lambda_0 ) = s_+ ( \lambda_0 ) - s_- (\lambda_0 )  $ for $s_{\pm} (\lambda_0 )  = \# \{ j \ ; \ \pm \mathrm{res} _{\lambda = \lambda_0 } \mu_j (\lambda ) >0 \} $.

\label{S3_lem_jump}
\end{lemma}

Proof.
In view of Lemma \ref{S3_lem_positiveB}, some eigenvalues $\mu_j (\lambda )$ have a pole at $\lambda_0 $.
If $ \pm \mathrm{res} _{\lambda = \lambda_0 } \mu_j (\lambda ) >0 $, we have $\mu_j (\lambda ) \to \mp \infty $ as $ \lambda \to \lambda_0 + 0 $ and $\mu_j (\lambda ) \to \pm \infty $ as $ \lambda \to \lambda_0 -0 $, respectively.
Then the number of negative eigenvalues decreases for $\mathrm{res} _{\lambda = \lambda_0 } \mu_j (\lambda ) <0 $ and increases for $\mathrm{res} _{\lambda = \lambda_0 } \mu_j (\lambda ) >0 $ when $\lambda $ passes through $\lambda_0 $ from $\alpha $.
This implies the lemma.
\qed

 \begin{lemma}
 
 If $\lambda_0 \in {\bf R} _{>0} $ is a pole of $\Lambda_k (\lambda )$, the residue $Q _{k,\mathcal{L} (\lambda_0 )} $ is negative.

 \label{S3_lem_residuepositive}
 \end{lemma}
 
 Proof. 
 Recall that $ B_k (\lambda_0 )$ is the subspace of $L^2 (\Gamma )$ spanned by $\partial _{\nu_k} \phi _{k,j} $ for $j\in \mathcal{L} (\lambda_0 )$.
 In view of (\ref{S2_eq_Pk}), we have for $0\not= f \in B_k (\lambda _0 ) $
 $$
 ( Q_{k,\mathcal{L} (\lambda_0 )} f ,f) _{L^2 (\Gamma )} = -\sum _{j\in \mathcal{L} (\lambda_0 )} | ( \partial _{\nu_k} \phi_{k,j} ,f ) _{L^2 (\Gamma )} | ^2 <0 .
 $$
Then we have $Q_{k,\mathcal{L} (\lambda_0 )} <0$. 
\qed  

\medskip

Let $ \lambda_0 \in \{ \lambda_{k,j} \} $.
We put $ m_k ( \lambda_0 ) = \mathrm{dim} \mathrm{Ran} Q_{k,\mathcal{L} (\lambda_0 )} $ and $m (\lambda_0 ) = \mathrm{dim} ( \mathrm{Ran} Q_{1,\mathcal{L} (\lambda_0)} \cap  \mathrm{Ran} Q_{2,\mathcal{L} (\lambda_0)}  )$.

 \begin{lemma}
 Let $\lambda _0 \in {\bf R} _{>0} $ be a pole of $B(\lambda )$. 
\\
 (1) If $\lambda_0 \not\in \{ \lambda_{1,j} \} \cap \{ \lambda_{2,j} \} $, we have $ \delta \mathcal{N}_{-\infty} (\lambda_0 ) = \gamma (m_2 (\lambda_0 )-m_1 (\lambda_0 ) ) $. \\
(2) If $\lambda _0 \in  \{ \lambda_{1,j} \} \cap \{ \lambda_{2,j} \} $, we have $ | \delta \mathcal{N}_{-\infty} (\lambda_0 ) - \gamma (m_2 (\lambda_0 ) -m_1 (\lambda_0 )) | \leq m (\lambda_0 ) $.
 
 \label{S3_lem_rankdeltan}
 \end{lemma}
 
 Proof.
 First we prove the assertion (1).
 Suppose $ \lambda_0 \in \{ \lambda _{1,j} \} $.
 We can expand $B(\lambda )$ in a small neighborhood of $\lambda_0 $ as 
 $$
 B(\lambda )= \frac{ \gamma D_{\Gamma} ^{(1+s)/4} Q_{1,\mathcal{L} (\lambda_0  )} D_{\Gamma} ^{(1+s)/4}  }{ \lambda _0 - \lambda  } + \widetilde{H} _{\lambda_0} (\lambda ) ,
 $$
 where $ \widetilde{H} _{\lambda_0} (\lambda ) $ is analytic. 
 From Lemma \ref{S3_lem_residuepositive}, we have $ Q_{1,\mathcal{L}  (\lambda_0 ) } < 0 $ and also $D_{\Gamma} ^{(1+s)/4} Q_{1,\mathcal{L} (\lambda_0  )} D_{\Gamma} ^{(1+s)/4} <0 $ so that $ D_{\Gamma} ^{(1+s)/4} Q_{1,\mathcal{L} (\lambda_0  )} D_{\Gamma} ^{(1+s)/4} $ has exactly $m_1 ( \lambda_0 )  $ strictly negative eigenvalues. 
Hence we have $\mathrm{sign} ( \mathrm{res} _{\lambda = \lambda_0 } \mu_j (\lambda ) )= - \gamma $.
In view of the assertion (2) in Lemma \ref{S3_lem_positiveB}, this means $s_+ (\lambda_0 )=0 $ and $ s _- (\lambda_0 )= m_1 (\lambda_0 )$ for $\gamma =1$, or $s_+ (\lambda_0 ) = m_1 (\lambda_0 )$ and $s_- (\lambda_0 )=0 $ for $\gamma =-1$.
Lemma \ref{S3_lem_jump} implies $ \delta \mathcal{N}_{-\infty} (\lambda_0 )= \gamma (m_2 (\lambda_0 )-m_1 (\lambda_0 ))$ with $ m_2 ( \lambda_0 ) =0 $. 
For the case $ \lambda_0 \in \{ \lambda_{2,j} \} $, we can see the same formula with $ m_1 ( \lambda_0 ) =0 $ by the similar way.
 Plugging these two cases, we obtain the assertion (1).

Let us prove the assertion (2). 
Suppose $\lambda _0 = \lambda_{1, i_1} = \lambda _{2,i_2}  $ for $\lambda_{1, i_1}  \in \{ \lambda_{1,j} \}  $ and $ \lambda _{2,i_2} \in   \{ \lambda_{2,j} \} $.
Then we have the following representation in a small neighborhood of $\lambda_0 $
$$
B(\lambda )= \frac{ \gamma Q _{\lambda_0}  }{\lambda_0 - \lambda } + \widetilde{H} _{\lambda_0 } (\lambda ) ,
$$
with $ Q_{\lambda_0} =  D_{\Gamma} ^{(1+s)/4} ( Q_{1,\mathcal{L} (\lambda_{1,i_1}  )} - Q_{2, \mathcal{L} (\lambda_{2,i_2} )}  ) D_{\Gamma} ^{(1+s)/4}$. 
We see that $ Q _{\lambda_0} <0 $ on $ B_1 ( \lambda _{1,i_1} ) \cap B_2 (\lambda _{2,i_2} )^{\perp } $ and $ Q _{\lambda_0} >0 $ on $ B_1 (\lambda _{1,i_1} )^{\perp} \cap B_2 (\lambda _{2,i_2} )  $.
If $\gamma =1$, we have $m_2 (\lambda_0 ) -m(\lambda_0 ) \leq s_+ (\lambda_0 ) \leq m_2 (\lambda_0 )$ and $m_1 (\lambda_0 ) -m(\lambda_0 )\leq s_- (\lambda_0 )\leq m_1 (\lambda_0 )$.
If $\gamma = -1$, we have $m_1 (\lambda_0 ) -m(\lambda_0 ) \leq s_+ (\lambda_0 ) \leq m_1 (\lambda_0 )$ and $m_2 (\lambda_0 ) -m(\lambda_0 )\leq s_- (\lambda_0 )\leq m_2 (\lambda_0 )$.
These inequalities and Lemma \ref{S3_lem_jump} imply the assertion (2).
\qed

\medskip

Now we have arrived at our main result on the Weyl type lower bound for $N_T (\lambda )$.

\begin{theorem}
We assume (A-1), one of (A-2-1) and (A-2-2), and (A-3).
For large $ \lambda \in {\bf  R}_{>0} $, we have 
\begin{equation}
N_T (\lambda ) \geq \gamma \sum_{\alpha < \lambda ' \leq \lambda} \left( m_1 ( \lambda ')  - m_2 ( \lambda ')  \right) -N_{-}(\alpha ) ,
\label{S3_eq_lower1}
\end{equation}
where the summation is taken over poles $ \lambda ' \in ( \alpha , \lambda ]$ of $\Lambda_1 (\lambda ) - \Lambda_2 (\lambda ) $.
Moreover, if $\gamma (V_1 -V_2 )>0$ where $V_1 , V_2 >0 $ are defined in (\ref{S3_eq_weylDirichlet}), $N_T (\lambda )$ satisfies asymptotically as $ \lambda \to \infty $ 
\begin{equation}
N_T ( \lambda ) \geq  \gamma ( V_1 - V_2 ) \lambda ^{d/2} +O(\lambda ^{(d-1)/2} ) .
\label{S3_eq_lower2} 
\end{equation}

\label{S3_thm_WeylNT}
\end{theorem}

Proof.
We prove for the case $\{ \lambda_{1,j} \} \cap \{ \lambda _{2,j} \} \not= \emptyset $.
For $\{ \lambda_{1,j} \} \cap \{ \lambda _{2,j} \} = \emptyset $, the proof is similar and can be slightly simplified.
Letting us recall that we call $\lambda$ is a singular ITE when $\lambda$ satisfies the latter condition of the assertion (2) in Lemma \ref{S2_lem_ITEeqcondition}, we put 
$$
N_{sng} (\lambda )= \# \{ \text{singular ITEs} \in ( \alpha , \lambda ] \subset {\bf R} _{>0} \} .
$$
Here $N_{sng} (\lambda )$ counts the number of singular ITEs with multiplicities taken into account. 
Note that $\mathcal{N} _{0}(\lambda)+N_{sng} (\lambda) \leq N_T (\lambda ) $ by the definition of $\mathcal{N}_0 (\lambda ) $ and Lemma \ref{S2_lem_ITEeqconditionB}. 
We take the summation of $ | \delta \mathcal{N}_{-\infty} (\lambda ' ) - \gamma (m_2 (\lambda ' ) -m_1 (\lambda ' )) | \leq m (\lambda ' ) $ in $ (\alpha , \lambda ]$.
Then we have 
$$
 \left| \mathcal{N}_{-\infty} (\lambda ) - \gamma \sum _{ \alpha < \lambda ' \leq \lambda } \left( m_2 ( \lambda ') - m_1  (\lambda ') \right)  \right| \leq N_{sng} (\lambda ) .
$$
See also Remark of Proposition \ref{S2_prop_green_meromor}.
Plugging this inequality and (\ref{S3_eq_diffNB}), we have 
\begin{gather*}
\begin{split}
&N_{-}(\lambda )-N_{-} (\alpha ) + \gamma \sum _{ \alpha < \lambda ' \leq \lambda } \left( m_1 ( \lambda ')  - m_2 ( \lambda ') \right) \\ 
&\leq \mathcal{N} _0 (\lambda ) + N_{sng} (\lambda )  \leq N_T (\lambda ) .
\end{split}
\end{gather*}
Since $N_{-} (\lambda ) \geq 0$, we obtain (\ref{S3_eq_lower1}). 

The inequality (\ref{S3_eq_lower1}) implies 
$$
N_T (\lambda ) \geq \gamma (N_1 (\lambda ) -N_2 (\lambda )) - N_- (\alpha ).
$$
The asymptotic estimate (\ref{S3_eq_lower2}) is a direct consequence of this inequality and Lemma \ref{S3_eq_weylDirichlet}.
\qed


\section{Non-scattering energy}

\subsection{Scattering theory for acoustic equations}
In the following, we derive a well-known scattering theory for the time-harmonic acoustic equation.
For the sake of simplicity, we consider the following operators :
$$
L = -n^{-1} \Delta , \quad L_0 = -\Delta \quad \text{on} \quad {\bf R}^d .
$$
Let $ \Omega = \mathrm{supp} (n(x)-1)$ be a bounded domain with smooth boundary.
We assume that $n\in C ({\bf R}^d)$, $n\big| _{\overline{\Omega}} \in C^{\infty} (\overline{\Omega} )$, and $n(x)$ is strictly positive for all $x\in {\bf R}^d$.
Moreover, we impose the following assumptions:

\medskip

{\bf (A)'} $n(x)=1$ and $ \partial _{\nu} n(x) \not= 0$ for all $x\in \partial \Omega $.

\medskip

The assumption (A)' corresponds (A-2) and (A-3) in this case.
The operators $L$ and $L_0$ are self-adjoint on $L^2 ({\bf R}^d , ndx )$ and $L^2 ({\bf R}^d ,dx)$, respectively, with the domain $ H^2 ({\bf R}^d)$.
For short, we use the notations $L^2_n ({\bf R}^d)= L^2 ({\bf R}^d , ndx)$ and $L^2 ({\bf R}^d)=L^2 ({\bf R}^d ,dx)$.
Obviously, we have $L \geq 0$ on $L^2_n ({\bf R}^d)$.

Let us list some basic facts which are well-known results in the spectral and the scattering theory.
For the Schr\"{o}dinger operators, see e.g. \cite{Ya} and \cite{Es}.
We can refer \cite{IsKu} and \cite{Mo} for the wave equations.
For the acoustic equation, the argument is similar.
We will omit the proofs.

\begin{lemma}
We have $\sigma_p (L)= \emptyset$ and $\sigma_{ac} (L)=\sigma_{ac} (L_0)= [0,\infty )$.

\label{S4_lem_essspec}
\end{lemma}

For the scattering theory, we consider the continuous spectrum.
Thus we take $\lambda >0$ in the following arguments.

Let $R(z)=(L-z)^{-1}$ and $R_0 (z)=(L_0 -z)^{-1}$ for $z\not\in [0,\infty )$.
We take a function $\chi \in C^{\infty} ({\bf R}^d)$ such that $ \chi (x)=1$ for $|x|>\rho +1$ and $\chi (x)=0$ for $|x|< \rho $ with a sufficiently large constant $\rho >1$.
In particular, we assume $\overline{\Omega} \subset \{ x\in {\bf R}^d \ ; \ |x| <\rho  \}$. 
Then we have 
\begin{gather}
\chi R(z)= R_0 (z) \chi - R_0 (z) (\chi L - L_0 \chi ) R(z) , \label{S4_eq_resolventeq1} \\
R(z) \chi = \chi R_0 (z) - R(z) (L \chi - \chi L_0 )R_0 (z) .\label{S4_eq_resolventeq2}
\end{gather}

In the following, $\mathcal{B}$ and $\mathcal{B}^*$ denote the pair of H\"{o}rmander's functional spaces (\cite{AgHo}).
In particular, the norm of $\mathcal{B}^*$ is given by
$$
\| u\| _{\mathcal{B}^*}^2 = \sup _{R>1} \frac{1}{R} \int_{|x|<R} |u(x)|^2 dx .
$$
Note that 
$$
\mathcal{B} \subset L^2 ({\bf R}^d ,dx) (\text{or} \ L^2 ({\bf R}^d ,ndx)) \subset \mathcal{B}^* .
$$
The space $ \mathcal{B}_0^*$ is the set of functions $u\in \mathcal{B}^*$ satisfying
$$
\lim_{R\to \infty} \frac{1}{R} \int _{|x|<R} |u(x)|^2 dx =0 .
$$

\begin{lemma}
For $ \lambda >0$, there exist the limits $ R( \lambda \pm i0 ) := \lim _{\epsilon \downarrow 0} R(\lambda \pm i\epsilon )$ in ${\bf B} (\mathcal{B} ; \mathcal{B}^* )$.
For any compact intervals $J\in (0,\infty)$, there exists a constant $C>0$ such that 
$$
\| R(\lambda \pm i0)f \| _{\mathcal{B}^*} \leq C \| f\| _{\mathcal{B}} ,
$$
for $ f\in \mathcal{B} $ where $\lambda $ varies on $J$.
Moreover, the mapping $ J\ni \lambda \mapsto (R(\lambda \pm i0)f,g)$ for $f,g \in \mathcal{B}$ is continuous.
$ R_0 (\lambda \pm i0 )$ satisfies the same kind of properties.
\label{S4_lem_LAP}
\end{lemma}

Note that $R_0 (\lambda \pm i0)$ is represented by the Green function :
$$
(R_0 (\lambda \pm i0 )f )(x)= \int _{{\bf R}^d} E(x-y; \lambda \pm i0) f(y) dy  , \quad f\in \mathcal{B} ,
$$
where $E(x;z) $ is given by
$$
E(x;z)= \frac{i}{4} \left( \frac{\sqrt{z}}{2\pi |x|} \right)^{(d-2)/2} h_{(d-2)/2}^{(1)} (\sqrt{z}|x|) .
$$
Here $h_{\alpha}^{(1)}  $ is the first Hankel function of order $ \alpha $ and the branch of $ \sqrt{z}$ is taken so that $ \mathrm{Im} \, \sqrt{z} >0$.

Let ${\bf h}_{\lambda} = L^2 ( S^{d-1} )$ with the inner product
$$
(\phi , \psi ) _{{\bf h}_{ \lambda }} = \frac{\lambda ^{(d-2)/2}}{2} \int _{S^{d-1}} \phi (\omega ) \overline{\psi (\omega )} d\omega . 
$$
Note that $L^2 ({\bf R}^d)$ is isometric to $\mathcal{H} := L^2 ((0,\infty ); {\bf h}_{\lambda} ; d\lambda )$.
We define the operator $ \mathcal{F}_0 (\lambda ) \in {\bf B} (\mathcal{B} ; {\bf h} _{\lambda} )$ by
$$
( \mathcal{F}_0 (\lambda )f)( \omega )= (2\pi )^{-d/2} \int _{{\bf R}^d} e^{-i\sqrt{\lambda} x\cdot \omega} f(x)dx , \quad \lambda >0 , \ \omega \in S^{d-1} ,
$$
for $ f\in \mathcal{B}$.
Thus we have
$$
( \mathcal{F}_0 (\lambda )^* \phi )(x)=  \frac{\lambda^{(d-2)/2}}{2^{(d+2)/2} \pi^{d/2}} \int _{S^{d-1}} e^{i\sqrt{\lambda} x\cdot \omega} \phi (\omega ) d\omega , \quad x\in {\bf R}^d ,
$$
for $ \phi \in {\bf h}_{\lambda} $.
Letting
$$
V=L\chi -\chi L_0,
$$
we define the distorted Fourier transformation by
$$
\mathcal{F}_{\pm} (\lambda )= \mathcal{F}_0 (\lambda ) \left( \chi-V^* R(\lambda \pm i0 ) \right) .
$$

\begin{lemma}
Let $\lambda >0$. \\
(1) We have $\mathcal{F}_{\pm} (\lambda) \in {\bf B} (\mathcal{B} ; {\bf h}_{\lambda} )$ and $\mathcal{F}_{\pm} (\lambda )^* \in {\bf B} ({\bf h}_{\lambda} ; \mathcal{B}^* )$.
\\
(2) We have $\mathcal{F}_{\pm} (\lambda ) \mathcal{B} = {\bf h}_{\lambda}$ and $\{ u\in \mathcal{B}^* \ ; \ (L-\lambda )u=0 \} = \mathcal{F}_{\pm} (\lambda )^* {\bf h} _{\lambda} $. \\
(3) For $f,g \in \mathcal{B}$, we have Stone's formula
$$
(R(\lambda + i0)f -R(\lambda -i0)f ,g)= 2\pi i (\mathcal{F}_{\pm} (\lambda )f , \mathcal{F}_{\pm} (\lambda )g)_{{\bf h} _{\lambda}} .
$$
(4) For $L_0$, $R_0 (\lambda \pm i0) $ and $ \mathcal{F}_0 (\lambda )$, the assertions (1)-(3) hold.
\label{S4_lem_distfourier}
\end{lemma}

Let $ u_{\pm}^{(0)} = R_0 (\lambda \pm i0)f$ and $u_{\pm} =R(\lambda \pm i0)f$ for $f\in  \mathcal{B}$.
These are unique solutions of the equations  $(L_0 -\lambda )u_{\pm}^{(0)} =f$ and $(L-\lambda )u_{\pm}=f$ with the Sommerfeld radiation condition 
$$ 
(\partial _r \mp i \sqrt{\lambda} )u^{(0)}_{\pm} ,(\partial _r \mp i \sqrt{\lambda })u_{\pm} \in \mathcal{B}_0^*,
$$
respectively.
Here $ \partial_r = \omega_x \cdot \nabla_x $ where $\omega_x = x/|x| \in S^{d-1} $.
Moreover, $ \mathcal{F}_0 (\lambda )$ and $\mathcal{F}_{\pm} (\lambda )$ appear in the far-field pattern of $u_{\pm}^{(0)}$ and $u_{\pm} $ in the sense of 
\begin{gather}
u_{\pm}^{(0)} (x)\sim C_{\pm} (\lambda ) |x|^{-(d-1)/2} e^{\pm i \sqrt{\lambda}|x|} (\mathcal{F}_0 (\lambda )f)( \pm \omega ) , \label{S4_eq_asymptotic1} \\
u_{\pm} (x)\sim C_{\pm} (\lambda ) |x|^{-(d-1)/2} e^{\pm i \sqrt{\lambda}|x|} (\mathcal{F}_{\pm} (\lambda )f)( \pm \omega ) ,\label{S4_eq_asymptotic2}
\end{gather}
as $|x|\to \infty $ in $\mathcal{B}_0^*$ where $ \omega = x/|x| \in S^{d-1}$ and $C_{\pm} (\lambda )= \pm \sqrt{\pi} \lambda ^{-1/4} e^{-i \pi (d-3)/4}$.
Then $ u_{\pm}^{(0)} $ and $u_{\pm} $ are outgoing for $+$ or incoming for $-$.

Let 
\begin{equation}
(\mathcal{F}_0 f)(\lambda )= \mathcal{F}_0 (\lambda )f, \quad (\mathcal{F}_{\pm} f)(\lambda )= \mathcal{F}_{\pm} (\lambda )f ,\quad f\in \mathcal{B}.
\label{S4_def_disFourier}
\end{equation}
Thus $ \mathcal{F}_0$ and $\mathcal{F}_{\pm} $ can be extended to a unitary operator from $L^2 ({\bf R}^d)$ or $L_n^2 ({\bf R}^d)$ to $ \mathcal{H}$.
The wave operators in view of the wave equation are defined by
\begin{equation}
W_{\pm} := {\mathop{{\rm s-lim}}_{t\to\pm \infty}} \, e^{it\sqrt{L}} \chi e^{-it \sqrt{L_0}} .
\label{S4_def_waveop}
\end{equation}
From the invariance property of the wave operators, $W_{\pm}$ can be represented by $ \mathcal{F}_{\pm} $ and $\mathcal{F}_0 $ as follows.

\begin{lemma}
The wave operators $W_{\pm}$ exist and complete i.e. $\mathrm{Ran} W_{\pm} = L^2 ({\bf R}^d)$.
Moreover, we have $ W_{\pm } = (\mathcal{F}_{\pm} )^* \mathcal{F}_0$.
\label{S4_lem_waveop}
\end{lemma}

The scattering operator is defined by $S=(W_+)^* W_-$.
We consider its Fourier transform $\widehat{S} = \mathcal{F}_0 S (\mathcal{F}_0)^*$.

\begin{lemma}
(1) We have a direct integral representation
$$
\widehat{S}= \int_0^{\infty} \oplus \widehat{S}(\lambda )d\lambda \quad \text{on} \quad \mathcal{H},
$$
where
\begin{equation}
\widehat{S} (\lambda )= 1-2\pi i A(\lambda ), \quad A(\lambda )= \mathcal{F}_+ (\lambda )V\mathcal{F}_0 (\lambda )^* .
\label{S4_eq_Smatrix}
\end{equation}
The S-matrix $\widehat{S} (\lambda )$ is unitary on ${\bf h} _{\lambda} $ for $\lambda >0$. \\
(2) For $\phi \in {\bf h}_{\lambda}$, we have 
$$
(\mathcal{F}_{-} (\lambda )^* \phi )(x) - ( \mathcal{F}_0 (\lambda )^* \phi )(x) \sim -C_{\pm} (\lambda ) |x|^{-(d-1)/2} e^{\pm i \sqrt{\lambda}|x|} (A (\lambda )\phi )( \omega ) ,
$$
as $|x|\to \infty $ in $\mathcal{B}_0 ^*$. 
\label{S4_lem_Smatrix}
\end{lemma}


\subsection{Layer potential}

In order to prove the equivalence of the S-matrix and the D-N map, we consider exterior and interior Dirichlet problems.
Thus we introduce Layer potentials. 
We follow the arguments in \cite{Es} or \cite{IsKu}.
We have to deal with Dirichlet eigenvalues, although we usually avoid them when we consider ISP.   
Then we slightly modify the arguments in view of the Laurent expansion of the D-N map as has been in \S 2.

We define the operators $ \delta , \delta _0 : L^2 (\partial \Omega )\to H^{-1/2} ({\bf R}^d)$ by 
$$
\int_{{\bf R}^d} (\delta f)(x) \overline{v(x)} n(x)dx = \int_{\partial \Omega} f(x')  \overline{v(x')} d\Sigma , 
$$
$$
\int_{{\bf R}^d} (\delta_0 f)(x) \overline{v(x)} dx = \int_{\partial \Omega} f(x') \overline{v(x')} d\Sigma ,
$$
for any $f\in L^2 (\partial \Omega)$ and $v\in H^{1/2} ({\bf R}^d)$, where $d\Sigma$ is the measure on $ \partial \Omega$.
Then $\delta^*, \delta_0^* : H^{1/2} ({\bf R}^d)\to L^2 (\partial \Omega)$ are the trace operators on $\partial \Omega$.
Since $ R(\lambda \pm i0)g \in H^2_{loc} ({\bf R}^d)$ for $g\in \mathcal{B} $, the mappings
$$
\mathcal{B} \ni g \mapsto \int _{\partial \Omega} f(x') \overline{(\delta^* R(\lambda \pm i0)g)(x') } d\Sigma ,
$$
$$
\mathcal{B} \ni g \mapsto \int _{\partial \Omega} f(x') \overline{(\delta^*_0 R_0 (\lambda \pm i0)g)(x') } d\Sigma , 
$$
for $ f\in L^2 (\partial \Omega )$ are bounded linear functionals.
Thus the operators $ R(\lambda \pm i0)\delta , R_0 (\lambda \pm i0)\delta_0  : L^2 (\partial \Omega) \to \mathcal{B}^*$ are defined by 
$$
\int _{{\bf R}^d} (R(\lambda \pm i0) \delta f)(x) \overline{g(x)} n(x)dx = \int _{\partial \Omega} f(x') \overline{ (\delta^* R(\lambda \mp i0)g)(x')} d\Sigma , 
$$
$$
\int _{{\bf R}^d} (R_0 (\lambda \pm i0) \delta_0 f)(x) \overline{g(x)} dx = \int _{\partial \Omega} f(x') \overline{ (\delta^*_0 R(\lambda \mp i0)g)(x')} d\Sigma , 
$$
for $g\in \mathcal{B}$.
Similarly, we define $R_0 (\lambda \pm i0)\delta_0 : L^2 ( \partial \Omega ) \to \mathcal{B}^* $.
Note that $R_0 (\lambda \pm i0)\delta$ is the well-known single layer potential :
$$
(R_0 (\lambda \pm i0)\delta_0 f)(x)= \int _{\partial \Omega} E(x-y'; \lambda \pm i0) f(y') d\Sigma (y').
$$
The integral on the right-hand side converges.
Hence $R_0 (\lambda \pm i0 )\delta_0 f$ is continuous for $f\in L^2 (\partial \Omega )$.
For a function $w(x)$, we put 
$$
w^+ (x)= \lim_{y\to x, y\in \Omega} w(y), \quad w^- (x)= \lim_{y\to x, y\in {\bf R}^d \setminus \Omega} w(y) , \quad x\in \partial \Omega .
$$
Then the jump relation on $\partial \Omega$ holds for $v= R_0 (\lambda \pm i0 )\delta_0  f$, $f\in L^2 (\partial \Omega)$ in the sense of 
$$
(\partial _{\nu} v)^+ -(\partial _{\nu} v)^- =f.
$$
The following jump relation of $R(\lambda \pm i0 )\delta$ also holds.

\begin{lemma}
Let $u_{\pm} = R(\lambda \pm i0 )\delta f$ for $f\in H^{3/2} (\partial \Omega)$.
Then we have 
$$
(\partial _{\nu} u_{\pm} )^+ - (\partial _{\nu} u_{\pm} )^- =f ,
$$
on $ \partial \Omega $.
\label{S4_lem_jump}
\end{lemma}

Proof. 
Note that $( -n^{-1} \Delta - \lambda ) u_{\pm} =0$ in ${\bf R}^d \setminus \partial \Omega $.
By the integration by parts, we have for any $v\in C_0^{\infty} ({\bf R}^d)$
$$
\int_{{\bf R}^d} u_{\pm} \cdot \overline{ (-n^{-1} \Delta -\lambda )v } ndx = \int _{\partial \Omega} \left( (\partial _{\nu} u_{\pm} )^+ -(\partial _{\nu} u_{\pm} )^- \right) \overline{v} d\Sigma .
$$
Putting $g= (-n^{-1}\Delta -\lambda )v$, we obtain 
$$
\int_{{\bf R}^d} u_{\pm} \cdot \overline{g } ndx = \int _{\partial \Omega} \left( (\partial _{\nu} u_{\pm} )^+ -(\partial _{\nu} u_{\pm} )^- \right) \overline{\delta^* R(\lambda \pm i0)g} d\Sigma .
$$
By the definition of $R(\lambda \pm i0)\delta$, we have $(\partial _{\nu} u_{\pm} )^+ - (\partial _{\nu} u_{\pm} )^- =f$.
\qed

\medskip

Now we introduce the exterior Dirichlet problem in $\Omega^e := {\bf R}^d \setminus \Omega $.
In the following, we use the notation $\mathcal{B}^* = \mathcal{B}^* (\Omega^e )$ which will not bring confusion.
Let $ L_e =-\Delta$ in $\Omega^e$ with the Dirichlet boundary condition on $\partial \Omega $.
For $ R_e (z)= (L_e -z)^{-1}$ for $z \not\in {\bf R}$, it is well-known the following facts.

\begin{lemma} 
For $ \lambda >0$, there exist the limits $ R_e ( \lambda \pm i0 ) := \lim _{\epsilon \downarrow 0} R_e (\lambda \pm i\epsilon )$ in ${\bf B} (\mathcal{B} ; \mathcal{B}^* )$.
For any compact intervals $J\in (0,\infty)$, there exists a constant $C>0$ such that 
$$
\| R_e (\lambda \pm i0)f \| _{\mathcal{B}^*} \leq C \| f\| _{\mathcal{B}} ,
$$
for $ f\in \mathcal{B}$ where $\lambda $ varies on $J$.
Moreover, the mapping $ J\ni \lambda \mapsto (R_e(\lambda \pm i0)f,g)$ for $f,g \in \mathcal{B}$ is continuous.

\label{S4_lem_extLAP}
\end{lemma}

Let $u^e_{\pm}\in \mathcal{B}^*$ be the outgoing (for $+$) or the incoming (for $-$) solution satisfying Sommerfeld's radiation condition of the equation
\begin{equation}
(-\Delta -\lambda )u^e_{\pm} =0 \quad \text{in} \quad \Omega^e , \quad u ^e _{\pm} \big| _{\partial \Omega} =f,
\label{S4_eq_extDirichlet}
\end{equation} 
with $ \lambda >0$.
The exterior D-N map $\Lambda_{\pm}^e (\lambda )$ is defined by 
\begin{equation}
\Lambda_{\pm}^e (\lambda )f= \partial _{\nu} u_{\pm}^e \quad \text{on} \quad \partial \Omega ,
\label{S4_def_extDN}
\end{equation}
where $\partial _{\nu} $ is the outward normal derivative on $ \partial \Omega $.
Note that $ u_{\pm}^e $ exist for $f\in H^{3/2} (\partial \Omega )$ as follows.
We can extend $f\in H^{3/2} (\partial \Omega)$ to $\widetilde{f} \in H^2 (\Omega^e )$ such that $\widetilde{f}\big|_{\partial \Omega} =f$ and $\widetilde{f} $ has a compact support.
Hence $u_{\pm}^e$ is given by
$$
u_{\pm}^e = \widetilde{f} - R_e (\lambda \pm i0)(-\Delta -\lambda )\widetilde{f} .
$$

The interior D-N map $\Lambda_n (\lambda )$ is defined by 
\begin{equation}
\Lambda_n (\lambda )f = \partial _{\nu} u^i \quad \text{on} \quad \partial \Omega ,
\label{S4_def_intDN}
\end{equation}
where $u^i$ is a unique solution (in a suitable subspace of $L^2 (\Omega )$) of the equation
\begin{equation}
(-n^{-1} \Delta -\lambda )u ^i =0 \quad \text{in} \quad \Omega , \quad u  ^i \big| _{\partial \Omega} =f.
\label{S4_eq_intDirichlet}
\end{equation} 
Replacing $ n(x)$ by $n_0 (x) :=1$, we also define the D-N maps 
$$
 \Lambda_{0} (\lambda )f= \partial _{\nu} u^i _0 ,  
$$
where $u_0^i $ is the solution of 
$$
(-\Delta -\lambda )u_0^i =0 \quad \text{in} \quad \Omega , \quad u_0^i \big| _{\partial \Omega} =f. 
$$

Let $ \sigma_D ( -n^{-1} \Delta )$ and $\sigma_D (-\Delta )$ be the sets of Dirichlet eigenvalues of $-n^{-1} \Delta$ and $-\Delta $ in $\Omega $.
As has been seen in \S 2.1, the D-N maps $\Lambda_n (\lambda )$ and $\Lambda_0 (\lambda ) $ have Laurent expansion in a small neighborhood of each Dirichlet eigenvalue of $-n^{-1} \Delta $ and $-\Delta $, respectively.
If $\lambda_0 \in \sigma_D (-n^{-1} \Delta )$ or $\lambda_0 \in \sigma_D (-\Delta )$, we denote by 
$$
\Lambda_{n} (\lambda )=\frac{Q_{\lambda_0}}{\lambda_0 -\lambda} + H_n (\lambda ) ,
$$
$$
\Lambda_{0} (\lambda )=\frac{Q_{0,\lambda_0}}{\lambda_0 -\lambda} + H_0 (\lambda ) ,
$$
the Laurent expansions of $\Lambda_n (\lambda )$ and $\Lambda_0 (\lambda )$ at $ \lambda_0$.

When $ \lambda_0\in \sigma_D (-n^{-1} \Delta )$, let $ E_n (\lambda _0 )$ be associated eigenspace of $-n^{-1}\Delta $.
The subspace $B_n (\lambda_0 ) $ of $L^2 (\partial \Omega )$ is spanned by $ \partial _{\nu} \phi _j \big| _{\partial \Omega} $ for all $\phi_j \in E_n (\lambda_0 ) $.
For $-\Delta$, we define $E_0 (\lambda_0 )$ and $B_0 (\lambda_0 )$ for $\lambda_0 \in \sigma_D (-\Delta )$ by the similar way.

In the following, we need to consider both of the cases where $\lambda$ is a Dirichlet eigenvalue or not.
Hence we define the following operators :
\begin{gather}
D_n (\lambda )= \left\{ 
\begin{split}
\Lambda_n (\lambda ) &, \quad \lambda \not\in \sigma_D ( -n^{-1} \Delta ) ,\\
H_n (\lambda ) &, \quad \lambda \in \sigma_D ( -n^{-1} \Delta) ,
\end{split} \right.
\label{S4_def_Dn}
\end{gather}
and 
\begin{gather}
D_0 (\lambda )= \left\{ 
\begin{split}
\Lambda_0 (\lambda ) &, \quad \lambda \not\in \sigma_D ( - \Delta) ,\\
H_0 (\lambda ) &, \quad \lambda \in \sigma_D ( - \Delta) .
\end{split} \right.
\label{S4_def_Dn}
\end{gather}
Then $ D_n (\lambda )   : H^{3/2} (\partial \Omega ) \to H^{1/2} (\partial \Omega )$ for $\lambda \not\in \sigma_D (-n^{-1} \Delta )$, and $ D_n (\lambda )   : H^{3/2} (\partial \Omega ) \cap B_n (\lambda ) ^c \to H^{1/2} (\partial \Omega ) \cap B_n (\lambda )^c$ for $ \lambda \in \sigma_D (-n^{-1} \Delta )$. 
Similarly, $ D_0 (\lambda )   : H^{3/2} (\partial \Omega ) \to H^{1/2} (\partial \Omega )$ for $\lambda \not\in \sigma_D (- \Delta )$, and $ D_0 (\lambda )   : H^{3/2} (\partial \Omega ) \cap B_0 (\lambda )^c \to H^{1/2} (\partial \Omega ) \cap B_0 (\lambda )^c $ for $ \lambda \in \sigma_D (- \Delta )$.

For $f\in H^{3/2} (\partial \Omega )$, we put 
$$
v_{\pm} = \chi^i u^i + \chi^e u_{\pm}^e,
$$
where $\chi^i$ and $\chi^e$ are the characteristic functions of $\Omega$ and $\Omega^e$, respectively.

\begin{lemma}
(1) Suppose $ \lambda >0$ is not a Dirichlet eigenvalue of $-n^{-1} \Delta$ in $\Omega$.
Then we have 
\begin{equation}
v_{\pm} = R(\lambda \pm i0) \delta (D_n (\lambda )-\Lambda_{\pm}^e (\lambda ) )f, 
\label{S4_eq_layerpotential_sol1}
\end{equation}
for $f\in H^{3/2} (\partial \Omega) $, $ \lambda \not\in \sigma_D (-n^{-1} \Delta )$ or $f\in H^{3/2} (\partial \Omega) \cap B_n (\lambda )^c $, $\lambda \in \sigma_D (-n^{-1} \Delta )$.\\
(2) We have 
\begin{gather*}
(D_n (\lambda )f,g)_{L^2 (\partial \Omega)} =(f,D_n (\lambda )g)_{L^2 (\partial \Omega)}  , \\
(\Lambda_{\pm}^e (\lambda )f,g)_{L^2 (\partial \Omega)} =(f,\Lambda_{\mp}^e (\lambda )g)_{L^2 (\partial \Omega)}  
\end{gather*}
for $f,g\in H^{3/2} (\partial \Omega) $, $ \lambda \not\in \sigma_D (-n^{-1} \Delta )$ or $f,g\in H^{3/2} (\partial \Omega) \cap B_n (\lambda )^c $, $\lambda \in \sigma_D (-n^{-1} \Delta )$. \\
(3) For $n(x)=n_0 (x)=1$, the assertions (1)-(3) hold, replacing $R(\lambda \pm i0)$, $D_n (\lambda )$ and $B_n (\lambda_0 )$ by $R_0 (\lambda \pm i0)$, $D_0 (\lambda )$ and $B_0 (\lambda_0 )$, respectively.
\label{S4_lem_layerpotential_sol}
\end{lemma}

Proof.
We shall show the lemma for $-n^{-1} \Delta $.
Suppose $\lambda \not\in \sigma_D (-n^{-1}\Delta )$.
Take an arbitrary function $v_0 \in C_0^{\infty} ({\bf R}^d )$.
By the integration by parts, we have 
\begin{gather*}
\begin{split}
\int _{B_{\rho}} v_{\pm} \cdot \overline{(-n^{-1} \Delta -\lambda )v_0 } ndx =& \, \int _{\partial \Omega} ( \partial _{\nu} u^i - \partial _{\nu} u_{\pm}^e ) \overline{v_0}  d\Sigma \\
& + \int _{S_{\rho}} (\partial _r u_{\pm}^e \cdot \overline{v_0} - u_{\pm}^e \cdot \overline{\partial_r v_0} ) dS_{\rho } ,
\end{split}
\end{gather*}
where  
$$
B_{\rho} = \{ x\in {\bf R}^d \ ; \ |x|<\rho \} , \quad S_{\rho} = \{ x\in {\bf R}^d \ ; \ |x|=\rho \} ,
$$
and $dS_{\rho} $ is the measure on $S_{\rho} $ induced from the Euclidean measure.
In view of the radiation condition, the second term on the right-hand side converges to zero as $\rho \to \infty $.
Then we have 
$$
\int _{{\bf R}^d} v_{\pm} \cdot \overline{(-n^{-1} \Delta -\lambda )v_0} n dx = \int _{\partial \Omega} ( \Lambda_n (\lambda )f- \Lambda_{\pm}^e (\lambda )f ) \overline{v_0} \, d\Sigma .
$$
Putting $g=(-n^{-1} \Delta -\lambda )v_0$, and using $v_0 =R(\lambda \mp i0)g$, we obtain (\ref{S4_eq_layerpotential_sol1}).
If $\lambda \in \sigma_D (-n^{-1}\Delta )$, we can obtain (\ref{S4_eq_layerpotential_sol1}) taking $f\in H^{3/2} (\partial \Omega ) \cap B_n (\lambda )^c$.

For the assertion (2), we consider the outgoing (for $+$) and the incoming (for $-$) solutions $v_+$ and $v_-$ of (\ref{S4_eq_extDirichlet}) with its boundary values $\delta ^* v_+  =f$ and $\delta ^* v_-  =g$. 
By the integration by parts, we have
\begin{gather*}
\begin{split}
&\int_{B_{\rho} \cap \Omega^e} \left(  (-n^{-1} \Delta -\lambda )v_+ \cdot \overline{v_-} - v_+ \cdot \overline{(-n^{-1}\Delta -\lambda )v_-} \right) n dx \\
&= \int _{\partial \Omega} \left( \Lambda_+^e (\lambda )f \cdot \overline{g} -f\cdot \overline{ \Lambda_-^e (\lambda )g } \right) d\Sigma \\
&\quad + \int _{S_{\rho} } \left( ( \partial_r -i\sqrt{\lambda} )v_+ \cdot \overline{v_-} - v_+ \cdot \overline{(\partial_r +i\sqrt{\lambda} ) v_-}   \right) dS_{\rho} .
\end{split}
\end{gather*}
Tending $\rho \to \infty$, we obtain the assertion (2) for $ \Lambda_{\pm}^e (\lambda )$.
For $D_n (\lambda )$, the proof is similar.
\qed


\subsection{Orthogonality of generalized eigenfunctions on the boundary}

\begin{lemma}
Let $\lambda \in \sigma_D (-n^{-1} \Delta )$.
Then $\delta ^* R(\lambda \pm i0)f \in B_n (\lambda )^c$ for $f\in \mathcal{B} $ if and only if $f\big| _{\Omega} \in E_n ( \lambda )^c $.

\label{S4_lem_orthogonalR}
\end{lemma}

Proof.
Note that $u_{\pm} = R(\lambda \pm i0)f$ satisfies 
$$
(-n^{-1}\Delta -\lambda ) u_{\pm} =f \quad \text{in} \quad \Omega , \quad \delta^* u_{\pm}  = \delta^* R(\lambda \pm i0)f .
$$
Take an arbitrary $v \in E_n (\lambda )$.
Then it follows from the integration by parts
$$
\int _{\Omega} f \cdot \overline{v} \, ndx = \int_{\partial \Omega} \delta ^* R(\lambda  \pm i0) f \cdot \overline{\partial _{\nu} v} \,d\Sigma.
$$
This equality implies the lemma.
\qed

\begin{lemma}
(1) Let $\lambda \in \sigma_D (-n^{-1}\Delta ) $.
Then we have $\delta ^* \mathcal{F}_{\pm} (\lambda )^* \phi  \in B_n (\lambda )^c$ for any $\phi \in {\bf  h}_{\lambda}$.
\\
(2) Let $\lambda \in \sigma_D (-\Delta ) $.
Then we have $\delta ^* _0 \mathcal{F}_{0} (\lambda )^* \phi  \in B_0 (\lambda )^c$ for any $\phi \in {\bf  h}_{\lambda}$.
\label{S4_lem_F0Fpm}
\end{lemma}

Proof.
Let $ u_{\pm} = \mathcal{F}_{\pm} (\lambda  )^* \phi $ for any $\phi \in {\bf h}_{\lambda} $.
It follows from the integration by parts in $\Omega$ that 
$$
\int _{\partial \Omega} \delta^* u_{\pm} \cdot \overline{\partial _{\nu} v} \, d\Sigma =0,
$$
for any $v\in E_n (\lambda )$.
Here we have used the equation $ (-n^{-1} \Delta -\lambda )u_{\pm} =0$.
Then we obtain the assertion (1).
For the assertion (2), the proof is similar.
\qed

\medskip

Let us introduce the operators $M_{\pm} (\lambda )$ and $ M_{\pm ,0} (\lambda )$ by 
\begin{gather*}
M_{\pm} (\lambda )f= \delta ^* R(\lambda \pm i0)\delta f , \\
M_{0,\pm} (\lambda )f= \delta ^*_0  R_0 (\lambda \pm i0)\delta_0 f ,
\end{gather*} 
for $f\in H^{1/2} (\partial \Omega )$.

\begin{lemma}
(1) Let $\lambda \not\in \sigma_D (-n^{-1} \Delta )$.
Then $M_{\pm} (\lambda ): H^{1/2} (\partial \Omega )\to H^{3/2} (\partial \Omega )$ is one to one. \\
(2) Let $\lambda \in \sigma_D (-n^{-1}\Delta )$.
Then $M_{\pm} (\lambda )$ is one to one as a mapping $H^{1/2} (\partial \Omega ) \cap B_n (\lambda) ^c \to H^{3/2} (\partial \Omega)\cap B_n (\lambda )^c$. \\
(3) Let $\lambda \not\in \sigma_D (-\Delta )$.
Then $M_{0,\pm} (\lambda ): H^{1/2} (\partial \Omega )\to H^{3/2} (\partial \Omega )$ is one to one. \\
(4) Let $\lambda\in \sigma_D (-\Delta )$.
Then $M_{0,\pm} (\lambda )$ is one to one as a mapping $H^{1/2} (\partial \Omega ) \cap B_0 (\lambda ) ^c \to H^{3/2} (\partial \Omega)\cap B_0 (\lambda )^c$.
\label{S4_lem_M}
\end{lemma}

Proof.
We shall prove (1) and (2).
For (3) and (4), the proof is similar.
Suppose $\lambda \not\in \sigma_D (-n^{-1} \Delta )$ and $M_{\pm} (\lambda ) f=0$.
Then $u_{\pm} = R(\lambda \pm i0) \delta f$ satisfies 
\begin{gather*}
(-\Delta - \lambda )u_{\pm} =0 \quad \text{in} \quad \Omega^e , \\
(-n^{-1} \Delta -\lambda )u_{\pm} =0 \quad \text{in} \quad \Omega , 
\end{gather*}
with the boundary condition $u_{\pm} \big| _{\partial \Omega} =0$.
Since $u_{\pm}$ is outgoing (for $+$) or incoming (for $-$), we have $u_{\pm} =0$ in $\Omega ^e $.
Moreover, $u_{\pm} $ is a Dirichlet eigenfunction of $-n^{-1} \Delta$. 
Then the assumption implies $ u_{\pm } =0$ in $\Omega$.
In view of Lemma \ref{S4_lem_jump}, we have $ (\partial _{\nu} u_{\pm} )^+ -(\partial _{\nu} u_{\pm} )^- = f =0$.

When $\lambda\in \sigma_D (-n^{-1}\Delta ) $, we can see $u_{\pm} =0$ in $\Omega ^e $ by the same way.
In $\Omega$, $ u_{\pm} $ is a Dirichlet eigenfunction and Lemma \ref{S4_lem_jump} implies $ (\partial _{\nu} u_{\pm})^+ =f$.
If $0\not= f\in B_n (\lambda )^c$, this is a contradiction.
Thus we have $f=0$ in $B_n (\lambda )^c $.
\qed

\begin{lemma}
(1) If $\lambda \not\in \sigma_D (-n^{-1}\Delta )$, $ D_n (\lambda ) - \Lambda_{\pm}^e (\lambda )  $ is an isomorphism from $ H^{3/2} (\partial \Omega )$ to $ H^{1/2} (\partial \Omega )$. 
If $\lambda \in \sigma_D (-n^{-1}\Delta )$, $  D_n (\lambda ) - \Lambda_{\pm}^e (\lambda ) $ is an isomorphism from the subspace $ H^{3/2} (\partial \Omega ) \cap B_n (\lambda )^c$ to the subspace $H^{1/2} (\partial \Omega ) \cap B_n (\lambda )^c$. \\
(2) If $\lambda \not\in \sigma_D (-\Delta )$, $ D_0 (\lambda ) - \Lambda_{\pm}^e (\lambda )  $ is an isomorphism from $ H^{3/2} (\partial \Omega )$ to $ H^{1/2} (\partial \Omega )$. 
If $\lambda \in \sigma_D (-\Delta )$, $  D_0 (\lambda ) - \Lambda_{\pm}^e (\lambda ) $ is an isomorphism from the subspace $ H^{3/2} (\partial \Omega ) \cap B_0 (\lambda )^c$ to the subspace $H^{1/2} (\partial \Omega ) \cap B_0 (\lambda )^c$. 
\label{S4_lem_MLambda}
\end{lemma}

Proof.
Let $\lambda \not\in \sigma_D (-n^{-1} \Delta )$.
It follows from Lemmas \ref{S4_lem_layerpotential_sol} and \ref{S4_lem_M} that 
$$
 M_{\pm} (\lambda )(D_n (\lambda )-\Lambda_{\pm}^e (\lambda ))=1 .
$$
Thus $M_{\pm} (\lambda ) : H^{1/2} (\partial \Omega ) \to H^{3/2} (\partial \Omega )$ is one to one and surjective.
In particular, $M_{\pm} (\lambda )$ is an isomorphism.
This shows the assertion (1) with $\lambda \not\in \sigma_D (-n^{-1} \Delta )$.
For the other cases, the proofs are completely parallel.
\qed


\subsection{From D-N map to S-matrix}

Let us derive two kinds of resolvent equations for $R_e (z)$.

\begin{lemma}
We have 
\begin{gather}
\chi R_e (\lambda \pm i0) = R_0 (\lambda \pm i0) \chi - R_0 (\lambda \pm i0) ( \chi L_e - L_0 \chi ) R_e (\lambda \pm i0 ),
\label{S4_eq_ext_resolventeq1prime} \\
R_e (\lambda \pm i0 ) \chi = \chi R_0 (\lambda \pm i0 ) - R_e (\lambda \pm i0 )(L_e \chi - \chi L_0 ) R_0 (\lambda \pm i0 ).\label{S4_eq_ext_resolventeq2prime}
\end{gather}
\label{S4_lem_ext_resolventeq_prime}
\end{lemma}

Proof.
The equation (\ref{S4_eq_ext_resolventeq1prime}) is a consequence of the equality
\begin{gather*}
\begin{split}
(L_0 -\lambda )\chi R_e (\lambda \pm i0 ) &= \chi (L_e -\lambda ) R_e (\lambda \pm i0 ) - (\chi L_e - L_0 \chi ) R_e (\lambda \pm i0 ) \\
&= \chi - (\chi L_e - L_0 \chi ) R_e (\lambda \pm i0 ).
\end{split}
\end{gather*}
Taking the adjoint, we also have (\ref{S4_eq_ext_resolventeq2}).
\qed

\begin{lemma}
Let $f\in \mathcal{B}  $. 
If $\lambda \not\in \sigma_D (-n^{-1} \Delta )$, we have 
\begin{equation}
R_e (\lambda \pm i0)(\chi^e f ) = R_0 (\lambda \pm i0) f- R(\lambda \pm i0)\delta (D_n (\lambda )- \Lambda_{\pm}^e (\lambda ))\delta^*_0 R_0 (\lambda \pm i0)f ,
\label{S4_eq_ext_resolventeq}
\end{equation}
and
\begin{equation}
R_e (\lambda \pm i0)f= R_0 (\lambda \pm i0)\left( 1- \delta_0 (D_n (\lambda )- \Lambda_{\pm}^e (\lambda ))\delta^* R (\lambda \pm i0) \right) (\chi^e f ) .
\label{S4_eq_ext_resolventeq2}
\end{equation}
If $\lambda \in \sigma_D (-n^{-1}\Delta )$, the equality (\ref{S4_eq_ext_resolventeq2}) holds for $f\in \mathcal{B}$ such that $f\big| _{\Omega} \in E_n (\lambda )^c $.

\label{S4_lem_extLAP}
\end{lemma}

Proof.
Let $v^e _{\pm} = R_e (\lambda \pm i0)f$ be the outgoing or incoming solution of 
$$
(-\Delta -\lambda )v^e _{\pm} =f \quad \text{in} \quad \Omega^e , \quad v^e _{\pm} \big| _{\partial \Omega} =0,
$$
for $f\in \mathcal{B} $.
For $\lambda \not\in \sigma_D (-n^{-1} \Delta )$, $w_{\pm} =R(\lambda \pm i0)\delta (D_n (\lambda )-\Lambda_{\pm}^e (\lambda ))g$ and $\widetilde{w}_{\pm} = R_0 (\lambda \pm i0 )f$ for $g \in H^{3/2} (\partial \Omega )$ satisfy 
$$
(-\Delta -\lambda )w_{\pm} =0 \quad \text{in} \quad \Omega^e , \quad w_{\pm} \big| _{\partial \Omega} =  g,
$$
and
$$
(-\Delta -\lambda ) \widetilde{w}_{\pm}= f \quad \text{in} \quad \Omega^e , \quad \widetilde{w}_{\pm} \big| _{\partial \Omega} =\delta_0^* R_0 (\lambda \pm i0 )f .
$$
Letting $g= \delta_0^* R_0 (\lambda \pm i0)f$, we obtain (\ref{S4_eq_ext_resolventeq}), since $v_{\pm}^e $, $w_{\pm}$ and $\widetilde{w}_{\pm} $ are outgoing (for $+$) or incoming (for $-$).
Taking the adjoint, we also have (\ref{S4_eq_ext_resolventeq2}).

Let us turn to the case $\lambda \in \sigma_D (-n^{-1}\Delta )$.
Take $0< \mu\not= \lambda$ in a sufficiently small neighborhood of $\lambda$.
Then we have $ \mu \not\in \sigma_D ( -n^{-1}\Delta )$ and  (\ref{S4_eq_ext_resolventeq2}) holds at $\mu $.
If we take $f\in \mathcal{B}$ such that $f\big| _{\Omega} \in E_n (\lambda )^c$, (\ref{S4_eq_ext_resolventeq2}) can be rewritten by 
$$
R_e (\mu \pm i0)f= R_0 (\mu \pm i0) f- R_0 (\mu \pm i0)\delta_0 (D_n (\mu )- \Lambda_{\pm}^e (\mu ))\delta^* R (\mu \pm i0)f ,
$$
from Lemma \ref{S4_lem_orthogonalR}.
Since $R_e (\mu \pm i0)$, $R_0 (\mu \pm i0)$ and $R(\mu \pm i0)$ are continuous in the weak $*$ sense in a neighborhood of $\lambda$, $\mu$ in the above equality can tend to $ \lambda$.
Thus we obtain (\ref{S4_eq_ext_resolventeq2}) at $\mu = \lambda $.
\qed

\medskip

We define 
\begin{equation}
\mathcal{F}_{\pm}^e (\lambda )= \mathcal{F}_0 (\lambda ) \left( \chi - ( \chi L_e - L_0 \chi ) R_e (\lambda \pm i0) \right) .
\label{S4_def_Fe+-}
\end{equation}
By the definition, we have $ \mathcal{F}_{\pm}^e (\lambda ) \in {\bf B} (\mathcal{B} ; {\bf h}_{\lambda} )$.

\begin{lemma}
We take a function $\phi \in {\bf h}_{\lambda} $.
Then $\mathcal{F}_{-} ^e (\lambda )^* \phi \in \mathcal{B}^*$ satisfies
$$
(-\Delta - \lambda ) \mathcal{F}_-^e (\lambda )^* \phi =0 \quad \text{in} \quad \Omega^e , \quad \mathcal{F}_-^e (\lambda )^* \phi  \big| _{\partial \Omega} =0.
$$
Moreover, $ \mathcal{F}_-^e (\lambda )^* \phi - \chi \mathcal{F}_0 (\lambda )^* \phi $ is outgoing and satisfies 
$$
\mathcal{F}_-^e (\lambda )^* \phi - \chi \mathcal{F}_0 (\lambda )^* \phi \sim -C_+ (\lambda ) |x|^{-(d-1)/2} e^{i\sqrt{\lambda}|x|}( A^e (\lambda )\phi )(\omega ) ,
$$
as $|x| \to \infty $ in $\mathcal{B}_0 ^*$ where $A^e (\lambda )= \mathcal{F}_+^e (\lambda ) (L_e \chi  - \chi L_0 ) \mathcal{F}_0 (\lambda )^*$. 
\label{S4_lem_exteriorEF}
\end{lemma}

Proof.
Note that $ \mathcal{F}_{-}^e (\lambda )^* \phi$ satisfies
\begin{equation}
\mathcal{F}_{-}^e (\lambda )^* \phi = \chi \mathcal{F}_0 (\lambda )^* \phi - R_e (\lambda + i0) ( L_e \chi - \chi L_0 ) \mathcal{F}_0 (\lambda )^* \phi .
\label{S4_def_Fe+-*}
\end{equation}
Thus (\ref{S4_eq_ext_resolventeq1prime}) and (\ref{S4_def_Fe+-*}) show 
$$
\mathcal{F}_-^e (\lambda )^* \phi - \chi \mathcal{F}_0 (\lambda )^* \phi \sim -C_+ (\lambda )|x|^{-(d-1)/2} e^{i\sqrt{\lambda}|x|} (A^e (\lambda )\phi )(\omega ),
$$
as $|x|\to \infty $ in $\mathcal{B}_0 ^*$.
\qed

\medskip

Now let us define the operators $\Gamma_{\pm} (\lambda ): H^{3/2} (\partial \Omega ) \to {\bf h}_{\lambda} $ by
\begin{equation}
\Gamma _{\pm} (\lambda )f = \mathcal{F}_0 (\lambda ) ( (-\Delta -\lambda ) ( \chi u_{\pm}^e )), 
\label{S4_def_gamma+-}
\end{equation}
where $ u_{\pm}^e \in \mathcal{B}^*$ is the outgoing (for $+$) or incoming (for $-$) solution of (\ref{S4_eq_extDirichlet}).
Obviously, $\Gamma _{\pm} (\lambda )$ depends only on $\Omega$.

\begin{lemma}
Let $u_{\pm}^e \in \mathcal{B}^* $ be the outgoing (for $+$) or incoming (for $-$) solution of (\ref{S4_eq_extDirichlet}).
We have for any $f\in H^{3/2} (\partial \Omega )$
\begin{equation}
u_{\pm}^e (x) \sim C_{\pm} (\lambda ) |x|^{-(d-1)/2} e^{\pm i \sqrt{\lambda}} (\Gamma _{\pm} (\lambda )f )(\pm \omega ),
\label{S4_eq_asymptotic_extGamma}
\end{equation}
as $|x|\to \infty $ in $ \mathcal{B}_0^*$. 
For $f\in H^{3/2} (\partial \Omega)$ with $\lambda \not\in \sigma_D  (-n^{-1}\Delta )$ or $f\in H^{3/2} (\partial \Omega) \cap B_n (\lambda )^c$ with $\lambda \in \sigma_D (-n^{-1}\Delta )$, $\Gamma _{\pm} (\lambda )$ is represented by
\begin{equation}
\Gamma _{\pm} (\lambda )f = \mathcal{F}_{\pm} (\lambda ) \delta (D_n (\lambda )-\Lambda  _{\pm}^e (\lambda )) f .
\label{S4_eq_Gammapm_formula}
\end{equation}

\label{S4_lem_Gammapm}
\end{lemma}

Proof.
In view of the equality
$$
(-\Delta -\lambda )\chi u_{\pm}^e = -2 \nabla \chi \cdot \nabla u_{\pm}^e - (\Delta \chi ) u_{\pm}^e =:g ,
$$
we have 
\begin{gather*}
\begin{split}
\chi (x)  u_{\pm}^e (x) &= (R_0 (\lambda \pm i0) g )(x) \\
&\sim C_{\pm} (\lambda ) |x|^{-(d-1)/2} e^{\pm i\sqrt{\lambda} |x|} (\mathcal{F}_0 (\lambda )g)(\pm \omega ),
\end{split}
\end{gather*}
as $|x| \to \infty$ in $\mathcal{B}_0^* $.
This shows (\ref{S4_eq_asymptotic_extGamma}).
Lemma \ref{S4_lem_layerpotential_sol} implies 
$$
u_{\pm}^e (x)\sim C_{\pm} (\lambda ) |x|^{-(d-1)/2} e^{\pm i\sqrt{\lambda} |x|} (\mathcal{F}_{\pm} (\lambda ) \delta (D_n (\lambda )-\Lambda  _{\pm}^e (\lambda )) f)(\pm \omega ),
$$
where $f$ is taken as in the lemma.
Since $ u_{\pm} ^e $ is outgoing or incoming, the uniqueness of the solution implies $\Gamma _{\pm} (\lambda ) f= \mathcal{F}_{\pm} (\lambda ) \delta (D_n (\lambda )-\Lambda  _{\pm}^e ( \lambda )) f $.
\qed

\medskip

\begin{lemma}
(1) $\Gamma _{\pm} (\lambda )$ is one to one on $H^{3/2} (\partial \Omega )$. \\
(2) The range of $ \Gamma _{\pm} (\lambda )^*  $ is dense in $L^2 (\partial \Omega )$.
\label{S4_lem_Gammapm_1to1}
\end{lemma}

Proof.
Suppose $ \Gamma _{\pm} (\lambda )f=0$ for some $f\in H^{3/2} (\partial \Omega )$.
In view of (\ref{S4_eq_asymptotic_extGamma}), we have $ u_{\pm}^e \sim 0$ in $\mathcal{B}_0^* $.
The Rellich's uniqueness theorem and the unique continuation property, we have $ u_{\pm}^e =0 $ in $\Omega^e $.
Thus we obtain $f=\delta^* u_{\pm}^e =0$.
This implies that $\Gamma _{\pm} (\lambda )$ is one to one.

Next we suppose $( \Gamma _{\pm} (\lambda )^* \phi , g)_{L^2 (\partial \Omega)} =0$ for any $ \phi \in {\bf h} _{\lambda} $.
The assertion (1) implies $g=0$.
Then we obtain the denseness of $\mathrm{Ran} \Gamma _{\pm} (\lambda )^*$ in $L^2 (\partial \Omega )$.
\qed

\begin{lemma}
We have $\Gamma_+ (\lambda )  M_+ (\lambda ) \Gamma_- (\lambda )^* = A^e (\lambda ) - A (\lambda )$.
In particular, $A(\lambda )$ and $ D_n (\lambda )$ determine each other.
\label{S4_lem_M+Ae}
\end{lemma}

Proof.
Let $ \lambda \not\in \sigma_D (-n^{-1} \Delta )$.
We put 
\begin{equation}
 u= \mathcal{F}_- (\lambda )^* \phi - \chi ^e  \mathcal{F}_-^e (\lambda )^* \phi ,
\label{S4_eq_DNtoS11}
\end{equation}
for any $\phi \in {\bf h} _{\lambda} $.
Thus $u$ satisfies 
$$
(-\Delta - \lambda )u=0 \quad \text{in} \quad \Omega ^e , \quad u\big| _{\partial \Omega} = \delta ^* \mathcal{F}_- (\lambda )^* \phi .
$$
Therefore, in view of Lemma \ref{S4_lem_layerpotential_sol}, $u$ can be represented by 
\begin{equation}
u=R(\lambda +i0 ) \delta ( D_n (\lambda ) - \Lambda_+^e (\lambda )) \delta^* \mathcal{F}_- (\lambda )^* \phi ,
\label{S4_eq_DNtoS12}
\end{equation}
in $\Omega^e $.
It follows from (\ref{S4_eq_DNtoS11}) that $u$ is outgoing and satisfies 
\begin{equation}
u(x) \sim C_+ (\lambda )|x|^{-(d-1)/2} e^{i\sqrt{\lambda}|x|} ( ( A^e (\lambda ) -  A(\lambda )) \phi )(\omega ) ,
\label{S4_eq_DNtoS13}
\end{equation}
as $|x|\to \infty$ in $\mathcal{B}_0^*$.
On the other hand, the representation (\ref{S4_eq_DNtoS12}) shows the asymptotic behavior 
\begin{equation}
u(x) \sim C_+ (\lambda )|x|^{-(d-1)/2} e^{i\sqrt{\lambda}|x|} (\mathcal{F}_+ (\lambda )  \delta (D_n (\lambda )-\Lambda_+ ^e (\lambda )) \delta^* \mathcal{F}_- (\lambda )^* \phi )(\omega ),
\label{S4_eq_DNtoS14}
\end{equation}
as $|x|\to \infty$ in $\mathcal{B}_0^*$.
Plugging (\ref{S4_eq_DNtoS13}) and (\ref{S4_eq_DNtoS14}), we obtain
$$
(A^e (\lambda )-A(\lambda ))\phi =\mathcal{F}_+ (\lambda )  \delta (D_n (\lambda )-\Lambda_+ ^e (\lambda )) \delta^* \mathcal{F}_- (\lambda )^* \phi .
$$
Inserting $ M_+ (\lambda )(D_n (\lambda )-\Lambda_+^e (\lambda ))=1$ on the right-hand side, we have 
$$
\mathcal{F}_+ (\lambda )  \delta (D_n (\lambda )-\Lambda_+ ^e (\lambda )) \delta^* \mathcal{F}_- (\lambda )^* \phi  = \Gamma_+ (\lambda ) M_+ (\lambda ) \Gamma_- (\lambda )^* \phi .
$$
We obtain 
\begin{equation}
A^e (\lambda )-A(\lambda )= \Gamma_+ (\lambda ) M_+ (\lambda ) \Gamma_- (\lambda ) .
\label{S4_eq_DNtoS15}
\end{equation}

Let us turn to the case $\lambda \in \sigma_D (-n^{-1} \Delta )$.
In view of Lemma \ref{S4_lem_F0Fpm}, we have $\delta ^* \mathcal{F}_- (\lambda )^* \phi \in B_n (\lambda )^c $ for any $\phi \in {\bf h}_{\lambda} $.
Then the operator
$$
\mathcal{F}_+ (\lambda )  \delta (D_n (\lambda )-\Lambda_+ ^e (\lambda )) \delta^* \mathcal{F}_- (\lambda )^* , 
$$
and the representation
$$
\Gamma_- (\lambda )^*  = (D_n (\lambda )-\Lambda_+^e (\lambda )) \delta ^* \mathcal{F}_- (\lambda )^* ,
$$
are well-defined on ${\bf h} _{\lambda} $.
Hence we can use Lemma \ref{S4_lem_MLambda} on $H^{3/2} (\partial \Omega)\cap B_n (\lambda )^c$ for this case.

Therefore, Lemma \ref{S4_lem_Gammapm_1to1} implies the lemma for any $\lambda >0$.
\qed

\subsection{Non-scattering energy}

Now we arrive at the main part on the non-scattering energies.
We consider the following boundary value problem :
\begin{gather}
(-n^{-1} \Delta -\lambda )u=0 \quad \text{in} \quad \Omega, \label{S4_eq_ITE1} \\
(- \Delta -\lambda )v=0 \quad \text{in} \quad \Omega, \label{S4_eq_ITE2} \\
u=v, \quad \partial _{\nu} u = \partial _{\nu} v \quad \text{on} \quad \partial \Omega . \label{S4_eq_ITE3}
\end{gather}
We denote by $\sigma_{NS} (L)$ the totality of NSEs of $L$. 
As has been mentioned in \S 1.2, Rellich's uniqueness theorem implies that $\lambda \in \sigma_{NS} (L)$ is an ITE associated with (\ref{S4_eq_ITE1})-(\ref{S4_eq_ITE3}).
We can also see the following fact.

\begin{lemma}
If $\lambda \in (\alpha , \infty )$ is a non-singular ITE associated with (\ref{S4_eq_ITE1})-(\ref{S4_eq_ITE3}), $\lambda $ is a NSE.

\label{S4_lem_NSEequivITE}
\end{lemma}

Proof.
Suppose that $\lambda \in (\alpha , \infty)$ is a non-singular ITE. 
Replacing $n(x)$ by $ n_0 (x) := 1$ on ${\bf R}^d$, we apply Lemma \ref{S4_lem_M+Ae} to $L_0$.
Since $ \Gamma _{\pm} (\lambda )$ and $A^e (\lambda )$ depend only on $\Omega$, we have the following formulas
\begin{gather*}
\Gamma _+ (\lambda ) M_{0,+} (\lambda ) \Gamma _- (\lambda )^* = A^e (\lambda ), \\
\Gamma _+ (\lambda ) M_{+} (\lambda ) \Gamma _- (\lambda )^* = A^e (\lambda )-A(\lambda ).
\end{gather*}
Talking the difference, we obtain 
\begin{equation}
\Gamma_+ (\lambda )( M_+ (\lambda ) - M_{0,+} (\lambda ) )\Gamma_- (\lambda ) ^* = -A(\lambda ) .
\label{S4_eq_diffDNtoS}
\end{equation}
Since $\Gamma _+ (\lambda )$ is one to one, we have $( M_+ (\lambda ) - M_{0,+} (\lambda ) )\Gamma_- (\lambda ) ^* \phi =0$ if and only if $A(\lambda )\phi=0$ for some $\phi \in {\bf h}_{\lambda} $.

Now we take $\lambda \in \sigma_{T,0} $.
Then there exists $0\not= f\in H^{3/2} (\partial \Omega )$ such that $ f\in \mathrm{Ker} ( D_n (\lambda ) -D_0 (\lambda ))$.
Putting $ g= (D_0 (\lambda )-\Lambda_+^e (\lambda ))f \in H^{1/2} (\partial \Omega )$, we have 
$$
(M_+ (\lambda ) -M_{0,+} (\lambda ))g = (D_n (\lambda )-\Lambda_+^e (\lambda ))^{-1} (D_0 (\lambda )-D_n (\lambda ))f=0.
$$
Note that $ \mathrm{Ker} ( D_n (\lambda ) - D_0 (\lambda ))$ is a subspace of $L^2 (\partial \Omega )$ with $ \mathrm{dim} \mathrm{Ker} ( D_n (\lambda ) - D_0 (\lambda )) \geq 1$.
In view of the assertion (2) of Lemma \ref{S4_lem_Gammapm_1to1},  there exists $0\not= \phi \in {\bf h} _{\lambda } $ such that $\Gamma_- (\lambda )^* \phi \in (D_0 (\lambda )-\Lambda_+^e (\lambda ) ) \mathrm{Ker}   ( D_n (\lambda ) - D_0 (\lambda ))$. 
Thus we have $A(\lambda )\phi =0$ which shows $\lambda \in \sigma_{NS} (L)$.
\qed

\medskip

We put 
$$
\gamma = \mathrm{sgn} ( \partial _{\nu} n \big| _{\partial \Omega} ) .
$$
For each $x\in \Omega $, we define
$$
V_n = (2\pi )^{-d} \mathrm{vol} (B_d ) \int_{\Omega} \sqrt{n(x)} dx, \quad
V_0 = (2\pi )^{-d}  \mathrm{vol} ( B_d ) \mathrm{vol} (\Omega ),
$$
where $B_d$ is the unit ball in ${\bf R}^d$.

\begin{theorem}
Let $\alpha >0$ be sufficiently small.
Suppose that the number of singular ITEs in $ ( \alpha , \lambda ] $ with multiplicities satisfies $o(\lambda ^{d/2} ) $ as $\lambda \to \infty $. 
If $ \gamma (V_n - V_0 )>0 $, we have 
$$
\# ( \sigma_{NS} (L) \cap (\alpha , \lambda ]) \geq \gamma (V_n - V_0 ) \lambda ^{d/2} + o (\lambda ^{d/2} ) ,
$$
as $\lambda \to \infty $.
\label{S4_thm_NSE}
\end{theorem}

Proof.
Theorem \ref{S3_thm_WeylNT} and its proof show the inequality
$$
N_T (\lambda ) \geq \mathcal{N}_0 (\lambda )+ N_{sng} (\lambda ) \geq  \gamma (V_n - V_0 ) \lambda ^{d/2} + O (\lambda ^{(d-1)/2} ) ,
$$
as $ \lambda \to \infty $.
Under the assumption $ N_{sng} (\lambda ) =o (\lambda ^{d/2}) $ as $\lambda \to \infty $, we obtain 
$$
\# \{ \text{non-singular ITEs in } (\alpha , \lambda ] \} \geq \gamma (V_n - V_0 ) \lambda ^{d/2} + o (\lambda ^{d/2} ) ,
$$
as $ \lambda \to \infty $.
Thus Lemma \ref{S4_lem_NSEequivITE} implies the theorem.
\qed

\end{document}